\documentclass{article}
\usepackage{maple2e}
 
\DefineParaStyle{Author}
\DefineParaStyle{Heading 1}
\DefineParaStyle{Heading 2}
\DefineParaStyle{Maple Output}
\DefineParaStyle{Title}
\DefineCharStyle{2D Math}
\DefineCharStyle{2D Output}
\begin{document}
\begin{maplegroup}
\begin{Title}
ON THE COMPLETE SOLUTION TO THE MOST GENERAL FIFTH DEGREE POLYNOMIAL
\end{Title}

\end{maplegroup}
\begin{maplegroup}
\begin{Author}
Richard J. Drociuk
\end{Author}

\end{maplegroup}
\begin{maplegroup}
\begin{center}
Physics Department
\end{center}

\end{maplegroup}
\begin{maplegroup}
\begin{center}
Simon Fraser University
\end{center}

\end{maplegroup}
\begin{maplegroup}
\begin{center}
Burnaby British Columbia, Canada.
\end{center}
\end{maplegroup}
\begin{maplegroup}
\begin{center}
April 10, 2000.
\end{center}
\end{maplegroup}
\begin{maplegroup}
\begin{center}
{\em Dedicated to Erland Samuel Bring \\ The first great pioneer into the
solution to the equation to
the fifth degree.}\end{center}

\end{maplegroup}
\begin{maplegroup}
\begin{Heading 1}
ABSTRACT   
\end{Heading 1}

\end{maplegroup}
\begin{maplegroup}
The motivation behind this note, is due to the non success in finding
the complete solution to the General Quintic Equation. The hope was to
have a solution with all the parameters precisely calculated in a
straight forward manner. This paper gives the closed form solution for
the five roots of the General Quintic Equation. They can be generated
on Maple V, or on the new version Maple VI. On the new
version of maple, Maple VI, it may be possible to insert all the
substitutions calculated in this paper, into one another, and
construct one large equation for the Tschirnhausian Transformation. The
solution also uses the Generalized Hypergeometric Function which Maple
V can calculate, robustly.

\end{maplegroup}
\begin{maplegroup}
\begin{Heading 1}
INTRODUCTION
\end{Heading 1}

\end{maplegroup}
\begin{maplegroup}
It has  been known since about 2000 BC, that the
Mesopotamians have been able to solve the Quadratic Equation with the
Quadratic Formula[Young, 1]. It took until 1545 AD, for Cardano to
publish his solution for the Cubic Equation, in his " Artis magnae
sive de regulis algebraicis". But it was actually Tartaglia who did
the original work to solve the cubic. Cardano's roommate, Ferrari (in
Cardano's Ars magna), solved the Quartic Equation at about the same
time Cardano solved the Cubic Equation. Tartaglia fought ferociously
against Cardano, Ferrari, and Sciopone Ferro, for stealing his
solution of the Cubic Equation. This situation was filled with
perjury, disputation, and bitterness. Finally, Cardano was thrown into
prison by the inquisition for heresy, for making the horoscope of
Christ[Guerlac, 2].

\end{maplegroup}
\begin{maplegroup}
Erland Samuel Bring (1786), was the first person to perform a
Tschirnhausian Transformation to a quintic equation, successfully. He
transformed a quintic with the fourth and third order terms missing,
i.e. x\symbol{94}5+px\symbol{94}2+qx+r=0, to the Bring Form
x\symbol{94}5-x-s=0 [Bring, 3]. This work was disputed by the
University of Lund, and was lost in the university's archives. I do
not know if an original copy still exists, there may still be one in
an observatory in Russia[Harley, 4]. It might be worth finding this
document, for history's sake, since I think Jerrard came along at a
later date, and claimed it as his own. The quest of a lot of the 19th
century mathematicians was to solve the Quintic Equation. Paolo
Ruffini (1803) gave a proof that the Quintic is not solvable with
radicals. Neils Henrik Abel(1824) gave a more rigorous proof, of the
same thing. Evartiste Galois(1830) invented group theory, and also
showed the same impossibility as Ruffini and Abel.  Group Theory
and Modular Functions would prove to be the mathematical framework by
which Bring's Equation was first solved [Young, 1], [Guerlac, 2]. In
1858, using Elliptic Modular Functions, Hermite solves Bring's
Equation. Kronecker, Gordan, Brioshi and Klein also gave solutions to
Bring's Equation closely after. For a good review and further
references see [Weisstein, 5], [King, 6] and [Klein, 7].

\end{maplegroup}
\begin{maplegroup}
 Of all the people that have solved Bring's Equation, or another
normal form, Klein's work seems to be the one that has the closest
thing to a complete solution to the General Quintic Equation. None of
the above solutions include Bring's Tranformation, to his normal form,
and leave too many parameters still to be calculated. I was looking
for a simple closed form solution which is easy to use, like the
Quadratic Formula, so I may substitute it into another set of
equations and formulas I am working on. This ruled out iteration
methods, and other approximation schemes. Then I looked at Modular
Function techniques, but these techniques leave too many parameters to
calculate, and are complicated with intricate steps that depend on the
properties of Modular Functions. Also, most solutions which use
Modular Functions, require a Modular Function still to be inverted
through the Hypergeometric Equation before a solution could be
obtained. Hermite's solution does not require this inversion. He does
calculate an elliptic nome, which he then inserts into a Modular Function
he defines. But it seems that these functions have a different period
that what he claimed. It also seems like the Russians or Weber realized
this and were doing same thing as Hermite with Weber's modular
functions, f1, f2 and f, but this also requires the inversion of f
[Prasolov, Solovyev, 8]. What is desirable, is to have just a two or
three step process to obtain the n roots of the nth degree polynomial
[Cockle, 9 and 10],[Harley, 11],[Cayley, 12]. So here we use only a three
step process to extract the roots from the General Quintic: 1) a
Tshirnhausian Transformation to Bring's equation; 2) the solution to a
generalized hypergeometric differential equation to solve Bring's
Equation; 3) and undo the Tschirnhausian Tranformation using Ferrari's
method to solve the Tschirnhausian Quartic. 

\end{maplegroup}
\begin{maplegroup}
\begin{Heading 1}
THE TSCHIRNHAUSIAN TRANSFORMATION TO BRING'S NORMAL FORM
\end{Heading 1}

\end{maplegroup}
\begin{maplegroup}
      The initial Tshirnhausian Transformation I use is a generalization
of Bring's [Bring, 3], but a simplification of Cayley's[Cayley, 13],
with a quartic substitution,

\end{maplegroup}
\begin{maplegroup}
\mapleresult
\begin{maplelatex}
\mapleinline{inert}{2d}{Tsh1 := x^4+d*x^3+c*x^2+b*x+a+y;}{%
\[
\mathit{Tsh1} := x^{4} + d\,x^{3} + c\,x^{2} + b\,x + a + y
\]
}
\end{maplelatex}

\end{maplegroup}
\begin{maplegroup}
to the General Quintic Equation,

\end{maplegroup}
\begin{maplegroup}
\mapleresult
\begin{maplelatex}
\mapleinline{inert}{2d}{Eq1 := x^5+m*x^4+n*x^3+p*x^2+q*x+r;}{%
\[
\mathit{Eq1} := x^{5} + m\,x^{4} + n\,x^{3} + p\,x^{2} + q\,x + r
\]
}
\end{maplelatex}

\end{maplegroup}
\begin{maplegroup}
Then by the process of elimination between Tsh1 and Eq1, the following
25 equations are obtained,

\end{maplegroup}
\begin{maplegroup}
\mapleresult
\begin{maplelatex}
\mapleinline{inert}{2d}{M15 := 1;}{%
\[
\mathit{M15} := 1
\]
}
\end{maplelatex}

\end{maplegroup}
\begin{maplegroup}
\mapleresult
\begin{maplelatex}
\mapleinline{inert}{2d}{M14 := d;}{%
\[
\mathit{M14} := d
\]
}
\end{maplelatex}

\end{maplegroup}
\begin{maplegroup}
\mapleresult
\begin{maplelatex}
\mapleinline{inert}{2d}{M13 := c;}{%
\[
\mathit{M13} := c
\]
}
\end{maplelatex}

\end{maplegroup}
\begin{maplegroup}
\mapleresult
\begin{maplelatex}
\mapleinline{inert}{2d}{M12 := b;}{%
\[
\mathit{M12} := b
\]
}
\end{maplelatex}

\end{maplegroup}
\begin{maplegroup}
\mapleresult
\begin{maplelatex}
\mapleinline{inert}{2d}{M11 := a+y;}{%
\[
\mathit{M11} := a + y
\]
}
\end{maplelatex}

\end{maplegroup}
\begin{maplegroup}
\mapleresult
\begin{maplelatex}
\mapleinline{inert}{2d}{M25 := m-d;}{%
\[
\mathit{M25} := m - d
\]
}
\end{maplelatex}

\end{maplegroup}
\begin{maplegroup}
\mapleresult
\begin{maplelatex}
\mapleinline{inert}{2d}{M24 := n-c;}{%
\[
\mathit{M24} := n - c
\]
}
\end{maplelatex}

\end{maplegroup}
\begin{maplegroup}
\mapleresult
\begin{maplelatex}
\mapleinline{inert}{2d}{M23 := p-b;}{%
\[
\mathit{M23} := p - b
\]
}
\end{maplelatex}

\end{maplegroup}
\begin{maplegroup}
\mapleresult
\begin{maplelatex}
\mapleinline{inert}{2d}{M22 := -y+q-a;}{%
\[
\mathit{M22} :=  - y + q - a
\]
}
\end{maplelatex}

\end{maplegroup}
\begin{maplegroup}
\mapleresult
\begin{maplelatex}
\mapleinline{inert}{2d}{M21 := r;}{%
\[
\mathit{M21} := r
\]
}
\end{maplelatex}

\end{maplegroup}
\begin{maplegroup}
\mapleresult
\begin{maplelatex}
\mapleinline{inert}{2d}{M35 := n+d*m-m^2-c;}{%
\[
\mathit{M35} := n + d\,m - m^{2} - c
\]
}
\end{maplelatex}

\end{maplegroup}
\begin{maplegroup}
\mapleresult
\begin{maplelatex}
\mapleinline{inert}{2d}{M34 := p-b-m*n+d*n;}{%
\[
\mathit{M34} := p - b - m\,n + d\,n
\]
}
\end{maplelatex}

\end{maplegroup}
\begin{maplegroup}
\mapleresult
\begin{maplelatex}
\mapleinline{inert}{2d}{M33 := q-a-m*p-y+d*p;}{%
\[
\mathit{M33} := q - a - m\,p - y + d\,p
\]
}
\end{maplelatex}

\end{maplegroup}
\begin{maplegroup}
\mapleresult
\begin{maplelatex}
\mapleinline{inert}{2d}{M32 := r-m*q+d*q;}{%
\[
\mathit{M32} := r - m\,q + d\,q
\]
}
\end{maplelatex}

\end{maplegroup}
\begin{maplegroup}
\mapleresult
\begin{maplelatex}
\mapleinline{inert}{2d}{M31 := d*r-m*r;}{%
\[
\mathit{M31} := d\,r - m\,r
\]
}
\end{maplelatex}

\end{maplegroup}
\begin{maplegroup}
\mapleresult
\begin{maplelatex}
\mapleinline{inert}{2d}{M45 := -c*m-m^3+b-p+d*m^2+2*m*n-d*n;}{%
\[
\mathit{M45} :=  - c\,m - m^{3} + b - p + d\,m^{2} + 2\,m\,n - d
\,n
\]
}
\end{maplelatex}

\end{maplegroup}
\begin{maplegroup}
\mapleresult
\begin{maplelatex}
\mapleinline{inert}{2d}{M44 := a+d*m*n+y-d*p-m^2*n-q+n^2+m*p-c*n;}{%
\[
\mathit{M44} := a + d\,m\,n + y - d\,p - m^{2}\,n - q + n^{2} + m
\,p - c\,n
\]
}
\end{maplelatex}

\end{maplegroup}
\begin{maplegroup}
\mapleresult
\begin{maplelatex}
\mapleinline{inert}{2d}{M43 := n*p-r+d*m*p-d*q+m*q-m^2*p-c*p;}{%
\[
\mathit{M43} := n\,p - r + d\,m\,p - d\,q + m\,q - m^{2}\,p - c\,
p
\]
}
\end{maplelatex}

\end{maplegroup}
\begin{maplegroup}
\mapleresult
\begin{maplelatex}
\mapleinline{inert}{2d}{M42 := m*r-d*r-m^2*q-c*q+d*m*q+n*q;}{%
\[
\mathit{M42} := m\,r - d\,r - m^{2}\,q - c\,q + d\,m\,q + n\,q
\]
}
\end{maplelatex}

\end{maplegroup}
\begin{maplegroup}
\mapleresult
\begin{maplelatex}
\mapleinline{inert}{2d}{M41 := -m^2*r-c*r+n*r+d*m*r;}{%
\[
\mathit{M41} :=  - m^{2}\,r - c\,r + n\,r + d\,m\,r
\]
}
\end{maplelatex}

\end{maplegroup}
\begin{maplegroup}
\mapleresult
\begin{maplelatex}
\mapleinline{inert}{2d}{M55 :=
b*m-2*m*p+q-y+c*n-2*d*m*n-a+3*m^2*n-c*m^2+d*m^3-n^2-m^4+d*p;}{%
\[
\mathit{M55} := b\,m - 2\,m\,p + q - y + c\,n - 2\,d\,m\,n - a + 
3\,m^{2}\,n - c\,m^{2} + d\,m^{3} - n^{2} - m^{4} + d\,p
\]
}
\end{maplelatex}

\end{maplegroup}
\begin{maplegroup}
\mapleresult
\begin{maplelatex}
\mapleinline{inert}{2d}{M54 :=
-d*m*p+d*m^2*n+c*p+d*q+2*m*n^2-c*m*n-m^3*n-2*n*p-m*q+m^2*p+b*n-d*n^2+r
;}{%
\maplemultiline{
\mathit{M54} :=  - d\,m\,p + d\,m^{2}\,n + c\,p + d\,q + 2\,m\,n
^{2} - c\,m\,n - m^{3}\,n - 2\,n\,p - m\,q + m^{2}\,p + b\,n \\
\mbox{} - d\,n^{2} + r }
}
\end{maplelatex}
\end{maplegroup}
\begin{maplegroup}
\mapleresult
\begin{maplelatex}
\mapleinline{inert}{2d}{M53 :=
c*q+2*m*n*p-d*n*p-p^2+b*p-n*q+d*m^2*p-m*r-d*m*q-c*m*p+m^2*q-m^3*p+d*r;
}{%
\maplemultiline{
\mathit{M53} := c\,q + 2\,m\,n\,p - d\,n\,p - p^{2} + b\,p - n\,q + d\,m^{2}\,p - m\,r - d\,m\,q - c\,m\,p + m^{2}\,q \\
\mbox{} - m^{3}\,p + d\,r }
}
\end{maplelatex}

\end{maplegroup}
\begin{maplegroup}
\mapleresult
\begin{maplelatex}
\mapleinline{inert}{2d}{M52 :=
b*q-c*m*q-n*r-d*m*r-d*n*q+c*r+m^2*r-m^3*q+2*m*n*q-p*q+d*m^2*q;}{%
\[
\mathit{M52} := b\,q - c\,m\,q - n\,r - d\,m\,r - d\,n\,q + c\,r
 + m^{2}\,r - m^{3}\,q + 2\,m\,n\,q - p\,q + d\,m^{2}\,q
\]
}
\end{maplelatex}

\end{maplegroup}
\begin{maplegroup}
\mapleresult
\begin{maplelatex}
\mapleinline{inert}{2d}{M51 := b*r-m^3*r-d*n*r+d*m^2*r+2*m*n*r-p*r-c*m*r;}{%
\[
\mathit{M51} := b\,r - m^{3}\,r - d\,n\,r + d\,m^{2}\,r + 2\,m\,n
\,r - p\,r - c\,m\,r
\]
}
\end{maplelatex}

\end{maplegroup}
\begin{maplegroup}
these equations are then substituted into the five by five matrix,

\end{maplegroup}
\begin{maplegroup}
\mapleresult
\begin{maplelatex}
\mapleinline{inert}{2d}{AA := matrix([[M11, M12, M13, M14, M15], [M21, M22, M23, M24, M25],
[M31, M32, M33, M34, M35], [M41, M42, M43, M44, M45], [M51, M52, M53,
M54, M55]]);}{%
\[
\mathit{AA} :=  \left[ 
{\begin{array}{ccccc}
\mathit{M11} & \mathit{M12} & \mathit{M13} & \mathit{M14} & 
\mathit{M15} \\
\mathit{M21} & \mathit{M22} & \mathit{M23} & \mathit{M24} & 
\mathit{M25} \\
\mathit{M31} & \mathit{M32} & \mathit{M33} & \mathit{M34} & 
\mathit{M35} \\
\mathit{M41} & \mathit{M42} & \mathit{M43} & \mathit{M44} & 
\mathit{M45} \\
\mathit{M51} & \mathit{M52} & \mathit{M53} & \mathit{M54} & 
\mathit{M55}
\end{array}}
 \right] 
\]
}
\end{maplelatex}

\end{maplegroup}
\begin{maplegroup}
taking the determinant of this matrix generates the polynomial which will
be reduced to Bring's Equation.
Let this polynomial be set to zero and solved, by first by transforming to
the Bring's Form. All the substitutions that are derived, transform poly into
Bring's form identically. Each step was checked and found to go to
zero, identically.

\end{maplegroup}

\begin{maplegroup}

 Let the transformed polynomial, poly = det(AA),  have the form,

\end{maplegroup}

\begin{maplegroup}
\mapleresult
\begin{maplelatex}
\mapleinline{inert}{2d}{poly := y^5+Poly4*y^4+Poly3*y^3+Poly2*y^2+A*y+B;}{%
\[
\mathit{poly} := y^{5} + \mathit{Poly4}\,y^{4} + \mathit{Poly3}\,
y^{3} + \mathit{Poly2}\,y^{2} + A\,y + B
\]
}
\end{maplelatex}

\end{maplegroup}
\begin{maplegroup}
setting Poly4 = 0, and solving for a,  gives

\end{maplegroup}
\begin{maplegroup}
\mapleresult
\begin{maplelatex}
\mapleinline{inert}{2d}{a :=
1/5*d*m^3+1/5*b*m-3/5*d*m*n-2/5*n^2+4/5*q-1/5*c*m^2+2/5*c*n-4/5*m*p+4/
5*m^2*n+3/5*d*p-1/5*m^4;}{%
\[
a := {\displaystyle \frac {1}{5}} \,d\,m^{3} + {\displaystyle 
\frac {1}{5}} \,b\,m - {\displaystyle \frac {3}{5}} \,d\,m\,n - 
{\displaystyle \frac {2}{5}} \,n^{2} + {\displaystyle \frac {4}{5
}} \,q - {\displaystyle \frac {1}{5}} \,c\,m^{2} + 
{\displaystyle \frac {2}{5}} \,c\,n - {\displaystyle \frac {4}{5}
} \,m\,p + {\displaystyle \frac {4}{5}} \,m^{2}\,n + 
{\displaystyle \frac {3}{5}} \,d\,p - {\displaystyle \frac {1}{5}
} \,m^{4}
\]
}
\end{maplelatex}

\end{maplegroup}
\begin{maplegroup}
Substituting a, and the substitutions,

\end{maplegroup}
\begin{maplegroup}
\mapleresult
\begin{maplelatex}
\mapleinline{inert}{2d}{b := alpha*d+xi;}{%
\[
b := \alpha \,d + \xi 
\]
}
\end{maplelatex}

\end{maplegroup}
\begin{maplegroup}
\mapleresult
\begin{maplelatex}
\mapleinline{inert}{2d}{c := d+eta;}{%
\[
c := d + \eta 
\]
}
\end{maplelatex}

\end{maplegroup}
\begin{maplegroup}
back into poly, and consider Poly3 and Poly4 as functions of d. Poly3
is quadratic in d, i.e.

\end{maplegroup}
\begin{maplegroup}
\mapleresult
\begin{maplelatex}
\mapleinline{inert}{2d}{Poly3 :=
((-2/5*m^2+n)*alpha^2+(17/5*m^2*n-17/5*m*p+4/5*m^3-2*n^2-13/5*n*m-4/5*
m^4+3*p+4*q)*alpha+22/5*m^2*p+21/5*m*p*n-19/5*p*n+12/5*n*m^4+5*r-2/5*m
^6-18/5*n^2*m^2+2*q+19/5*n^2*m-4*m^3*n+8/5*m^2*n+3*q*m^2-3*m*r-3/5*p^2
-3*n*q-5*q*m-2*m*p-2/5*m^4-12/5*m^3*p+n^3+4/5*m^5-3/5*n^2)*d^2+((-4/5*
m^2*xi+21/5*m^2*p-21/5*m^3*n-5*p*n-13/5*eta*m*n+2*n*xi+4/5*eta*m^3-21/
5*q*m+23/5*n^2*m+4/5*m^5+3*eta*p+5*r)*alpha+26/5*q*m^2-26/5*m^3*p+52/5
*m*p*n+4*q*xi-6*m*r-22/5*n*q-7*n^2*m^2-4/5*m^4*xi+17/5*m^2*n*xi-17/5*m
*p*xi-4/5*eta*m^4+4*eta*q-6/5*eta*n^2-4/5*m^6+7*m^2*r+4/5*m^7-3*p^2-2*
n^2*xi-31/5*q*m^3+28/5*m^4*p-4*eta*m*p+4/5*m^5*eta+19/5*n^2*m*eta-28/5
*m^5*n+5*r*eta+16/5*eta*m^2*n+58/5*q*m*n-81/5*m^2*p*n+56/5*n^2*m^3-23/
5*p*q+23/5*m*p^2-13/5*n*m*xi+3*p*xi+4/5*m^3*xi+22/5*p*m^2*eta-5*m*q*et
a-19/5*p*n*eta-29/5*n^3*m-7*n*r+24/5*n*m^4+29/5*p*n^2+6/5*n^3-4*m^3*n*
eta)*d+5*r*xi+n*xi^2+16/5*q*m^4+2*q*eta^2-5*p*n*xi-7*n^2*m^2*eta-21/5*
q*m*xi+52/5*m*p*n*eta-2/5*m^4*eta^2+21/5*m^2*p*xi+4/5*m^5*xi-2/5*m^2*x
i^2-2/5*q^2+16/5*m^6*n+26/5*q*m^2*eta-4*p*r-26/5*m^3*p*eta-22/5*q*n*et
a+6/5*n^3*eta-3*p^2*eta+32/5*n^3*m^2-2*m*p*eta^2-21/5*m^3*n*xi+4/5*m^3
*eta*xi-3/5*n^4-3/5*n^2*eta^2+23/5*n^2*m*xi+3*p*eta*xi+12/5*n^2*q-2/5*
m^8-52/5*m*p*n^2-13/5*n*m*eta*xi-4/5*m^6*eta-4*m^3*r-16/5*m^5*p-22/5*m
^2*p^2+24/5*n*eta*m^4+24/5*m*p*q+8/5*n*m^2*eta^2-44/5*m^2*n*q+8*m*n*r+
4*n*p^2-8*n^2*m^4+64/5*m^3*p*n-6*m*r*eta;}{%
\maplemultiline{
\mathit{Poly3} := (( - {\displaystyle \frac {2}{5}} \,m^{2} + n)
\,\alpha ^{2} + ({\displaystyle \frac {17}{5}} \,m^{2}\,n - 
{\displaystyle \frac {17}{5}} \,m\,p + {\displaystyle \frac {4}{5
}} \,m^{3} - 2\,n^{2} - {\displaystyle \frac {13}{5}} \,n\,m - 
{\displaystyle \frac {4}{5}} \,m^{4} + 3\,p + 4\,q)\,\alpha  \\
\mbox{} + {\displaystyle \frac {22}{5}} \,m^{2}\,p + 
{\displaystyle \frac {21}{5}} \,m\,p\,n - {\displaystyle \frac {
19}{5}} \,p\,n + {\displaystyle \frac {12}{5}} \,n\,m^{4} + 5\,r - {\displaystyle \frac {2}{5}} \,m^{6} - 
{\displaystyle \frac {18}{5}} \,n^{2}\,m^{2} + 2\,q + 
{\displaystyle \frac {19}{5}} \,n^{2}\,m \\
\mbox{} - 4\,m^{3}\,n + {\displaystyle \frac {8}{5}} \,m^{2}\,n + 3\,q\,m^{2} - 3\,m\,r - {\displaystyle \frac {3}{5}} \,p^{2} - 3\,n\,q - 5\,q\,m - 2\,m\,p - {\displaystyle \frac {2}{5}} \,m^{4} \\
\mbox{} - {\displaystyle \frac {12}{5}} \,m^{3}\,p + n^{3} + 
{\displaystyle \frac {4}{5}} \,m^{5} - {\displaystyle \frac {3}{5
}} \,n^{2})d^{2}\mbox{} + (( - {\displaystyle \frac {4}{5}} \,m^{
2}\,\xi  + {\displaystyle \frac {21}{5}} \,m^{2}\,p - {\displaystyle \frac {21}{5}} \,m^{3}\,n - 5\,p\,n \\
\mbox{} - {\displaystyle \frac {13}{5}} \,\eta \,m\,n + 2\,n\,\xi
  + {\displaystyle \frac {4}{5}} \,\eta \,m^{3} - {\displaystyle \frac {21}{5}} \,q\,m + {\displaystyle 
\frac {23}{5}} \,n^{2}\,m + {\displaystyle \frac {4}{5}} \,m^{5} + 3\,\eta \,p + 5\,r)\alpha \mbox{} + {\displaystyle 
\frac {26}{5}} \,q\,m^{2} \\
\mbox{} - {\displaystyle \frac {26}{5}} \,m^{3}\,p + {\displaystyle \frac {52}{5}} \,m\,p\,n + 4\,q\,\xi  - 6\,m\,r - {\displaystyle \frac {22}{5}} \,n\,q - 7\,n^{2}\,m^{2} - {\displaystyle \frac {4}{5}} \,m^{4}\,\xi 
 + {\displaystyle \frac {17}{5}} \,m^{2}\,n\,\xi  \\
\mbox{} - {\displaystyle \frac {17}{5}} \,m\,p\,\xi  - 
{\displaystyle \frac {4}{5}} \,\eta \,m^{4} + 4\,\eta \,q - 
{\displaystyle \frac {6}{5}} \,\eta \,n^{2} - {\displaystyle 
\frac {4}{5}} \,m^{6} + 7\,m^{2}\,r + {\displaystyle \frac {4}{5}
} \,m^{7} - 3\,p^{2} - 2\,n^{2}\,\xi  \\
\mbox{} - {\displaystyle \frac {31}{5}} \,q\,m^{3} + 
{\displaystyle \frac {28}{5}} \,m^{4}\,p - 4\,\eta \,m\,p + 
{\displaystyle \frac {4}{5}} \,m^{5}\,\eta  + {\displaystyle 
\frac {19}{5}} \,n^{2}\,m\,\eta  - {\displaystyle \frac {28}{5}} 
\,m^{5}\,n + 5\,r\,\eta  \\
\mbox{} + {\displaystyle \frac {16}{5}} \,\eta \,m^{2}\,n + 
{\displaystyle \frac {58}{5}} \,q\,m\,n - {\displaystyle \frac {
81}{5}} \,m^{2}\,p\,n + {\displaystyle \frac {56}{5}} \,n^{2}\,m^{3} - 
{\displaystyle \frac {23}{5}} \,p\,q + {\displaystyle \frac {23}{5}} \,m\,p^{2} - {\displaystyle \frac {13}{5}} \,n\,m\,\xi  \\
\mbox{} + 3\,p\,\xi  + {\displaystyle \frac {4}{5}} \,m^{3}\,\xi  + 
{\displaystyle \frac {22}{5}} \,p\,m^{2}\,\eta  - 5\,m\,q\,\eta  - {\displaystyle \frac {19}{5}} \,p\,n\,\eta  - {\displaystyle \frac {29}{5}} \,n^{3}\,m - 7\,n\,r + {\displaystyle \frac {24}{5}} \,n\,m^{4} \\
\mbox{} + {\displaystyle \frac {29}{5}} \,p\,n^{2} + {\displaystyle \frac {6}{5}} \,n^{3} - 4\,m^{3}\,n\,\eta )d\mbox{} + 5\,r\,\xi  + n\,\xi ^{2} + {\displaystyle \frac {16}{5}} \,q\,m^{4} + 2\,q\,\eta ^{2} - 5\,p\,n\,\xi  \\
\mbox{} - 7\,n^{2}\,m^{2}\,\eta  - {\displaystyle \frac {21}{5}} 
\,q\,m\,\xi  + {\displaystyle \frac {52}{5}} \,m\,p\,n\,\eta  - {\displaystyle \frac {2}{5}} \,m^{4}\,\eta ^{2} + 
{\displaystyle \frac {21}{5}} \,m^{2}\,p\,\xi  + {\displaystyle \frac {4}{5}} \,m^{5}\,\xi  - 
{\displaystyle \frac {2}{5}} \,m^{2}\,\xi ^{2} \\
\mbox{} - {\displaystyle \frac {2}{5}} \,q^{2} + {\displaystyle 
\frac {16}{5}} \,m^{6}\,n + {\displaystyle \frac {26}{5}} \,q\,m^{2}\,\eta  - 4\,p\,r - {\displaystyle \frac {26}{5}} \,m^{3}\,p\,\eta  - {\displaystyle \frac {22}{5}} \,q\,n\,\eta  + {\displaystyle \frac {6}{5}} \,n^{3}\,\eta  - 3\,p^{2}
\,\eta  \\
\mbox{} + {\displaystyle \frac {32}{5}} \,n^{3}\,m^{2} - 2\,m\,p\,\eta ^{2} - {\displaystyle \frac {21}{5}} \,m^{3}\,n\,\xi  + {\displaystyle \frac {4}{5}} \,m^{3}\,\eta \,\xi  - 
{\displaystyle \frac {3}{5}} \,n^{4} - {\displaystyle \frac {3}{5
}} \,n^{2}\,\eta ^{2} + {\displaystyle \frac {23}{5}} \,n^{2}\,m\,\xi  \\
\mbox{} + 3\,p\,\eta \,\xi  + {\displaystyle \frac {12}{5}} \,n^{2}\,q - {\displaystyle \frac {2}{5}} \,m^{8} - {\displaystyle \frac {
52}{5}} \,m\,p\,n^{2} - {\displaystyle \frac {13}{5}} \,n\,m\,\eta \,\xi  - {\displaystyle \frac {4}{5}} \,m^{6}\,\eta  - 4
\,m^{3}\,r - {\displaystyle \frac {16}{5}} \,m^{5}\,p \\
\mbox{} - {\displaystyle \frac {22}{5}} \,m^{2}\,p^{2} + 
{\displaystyle \frac {24}{5}} \,n\,\eta \,m^{4} + {\displaystyle \frac {24}{5}} \,m\,p\,q + {\displaystyle \frac {8}{5}} \,n\,m^{2}\,\eta ^{2} - {\displaystyle \frac {44}{5}} \,m^{2}\,n\,q + 8\,m\,n\,r + 4\,n\,p^{2} \\
\ }
}
\end{maplelatex}

\end{maplegroup}
\begin{maplegroup}
\mapleresult
\begin{maplelatex}
\mapleinline{inert}{2d}{Poly3 :=
((-2/5*m^2+n)*alpha^2+(17/5*m^2*n-17/5*m*p+4/5*m^3-2*n^2-13/5*n*m-4/5*
m^4+3*p+4*q)*alpha+22/5*m^2*p+21/5*m*p*n-19/5*p*n+12/5*n*m^4+5*r-2/5*m
^6-18/5*n^2*m^2+2*q+19/5*n^2*m-4*m^3*n+8/5*m^2*n+3*q*m^2-3*m*r-3/5*p^2
-3*n*q-5*q*m-2*m*p-2/5*m^4-12/5*m^3*p+n^3+4/5*m^5-3/5*n^2)*d^2+((-4/5*
m^2*xi+21/5*m^2*p-21/5*m^3*n-5*p*n-13/5*eta*m*n+2*n*xi+4/5*eta*m^3-21/
5*q*m+23/5*n^2*m+4/5*m^5+3*eta*p+5*r)*alpha+26/5*q*m^2-26/5*m^3*p+52/5
*m*p*n+4*q*xi-6*m*r-22/5*n*q-7*n^2*m^2-4/5*m^4*xi+17/5*m^2*n*xi-17/5*m
*p*xi-4/5*eta*m^4+4*eta*q-6/5*eta*n^2-4/5*m^6+7*m^2*r+4/5*m^7-3*p^2-2*
n^2*xi-31/5*q*m^3+28/5*m^4*p-4*eta*m*p+4/5*m^5*eta+19/5*n^2*m*eta-28/5
*m^5*n+5*r*eta+16/5*eta*m^2*n+58/5*q*m*n-81/5*m^2*p*n+56/5*n^2*m^3-23/
5*p*q+23/5*m*p^2-13/5*n*m*xi+3*p*xi+4/5*m^3*xi+22/5*p*m^2*eta-5*m*q*et
a-19/5*p*n*eta-29/5*n^3*m-7*n*r+24/5*n*m^4+29/5*p*n^2+6/5*n^3-4*m^3*n*
eta)*d+5*r*xi+n*xi^2+16/5*q*m^4+2*q*eta^2-5*p*n*xi-7*n^2*m^2*eta-21/5*
q*m*xi+52/5*m*p*n*eta-2/5*m^4*eta^2+21/5*m^2*p*xi+4/5*m^5*xi-2/5*m^2*x
i^2-2/5*q^2+16/5*m^6*n+26/5*q*m^2*eta-4*p*r-26/5*m^3*p*eta-22/5*q*n*et
a+6/5*n^3*eta-3*p^2*eta+32/5*n^3*m^2-2*m*p*eta^2-21/5*m^3*n*xi+4/5*m^3
*eta*xi-3/5*n^4-3/5*n^2*eta^2+23/5*n^2*m*xi+3*p*eta*xi+12/5*n^2*q-2/5*
m^8-52/5*m*p*n^2-13/5*n*m*eta*xi-4/5*m^6*eta-4*m^3*r-16/5*m^5*p-22/5*m
^2*p^2+24/5*n*eta*m^4+24/5*m*p*q+8/5*n*m^2*eta^2-44/5*m^2*n*q+8*m*n*r+
4*n*p^2-8*n^2*m^4+64/5*m^3*p*n-6*m*r*eta;}{%
\maplemultiline{
\mbox{} - 8\,n^{2}\,m^{4} + {\displaystyle \frac {64}{5}} \,m^{3}
\,p\,n - 6\,m\,r\,\eta  }
}
\end{maplelatex}

\end{maplegroup}
\begin{maplegroup}
Setting Poly3 to zero by setting each coefficient multiplying each
power of d equal to zero. The d\symbol{94}2 term is multiplied by a
Quadratic Equation in alpha, solving for alpha gives,

\end{maplegroup}
\begin{maplegroup}
\mapleresult
\begin{maplelatex}
\mapleinline{inert}{2d}{alpha :=
1/2*(-13*n*m-10*n^2+4*m^3+20*q+17*m^2*n-4*m^4+15*p-17*m*p+sqrt(-40*q*m
^4+80*q*m^2+40*m^3*p+60*n*p^2-15*n^2*m^4-190*m*p*n-200*n*q-15*n^2*m^2+
400*q^2+60*n^3*m^2-100*n^2*q-80*m*p*n^2+200*m^2*r+225*p^2-120*m^3*r+40
*m^5*p+265*m^2*p^2-40*q*m^3-80*m^4*p-20*q*m*n+360*m^2*p*n+30*n^2*m^3+6
00*p*q-510*m*p^2-120*n^3*m-680*m*p*q+260*m^2*n*q+300*m*n*r-500*n*r+80*
p*n^2-170*m^3*p*n+60*n^3))/(2*m^2-5*n);}{%
\maplemultiline{
\alpha  := {\displaystyle \frac {1}{2}} ( - 13\,n\,m - 10\,n^{2}
 + 4\,m^{3} + 20\,q + 17\,m^{2}\,n - 4\,m^{4} + 15\,p - 17\,m\,p
 + \mathrm{sqrt}( - 40\,q\,m^{4} \\
\mbox{} + 80\,q\,m^{2} + 40\,m^{3}\,p + 60\,n\,p^{2} - 15\,n^{2}
\,m^{4} - 190\,m\,p\,n - 200\,n\,q - 15\,n^{2}\,m^{2} \\
\mbox{} + 400\,q^{2} + 60\,n^{3}\,m^{2} - 100\,n^{2}\,q - 80\,m\,
p\,n^{2} + 200\,m^{2}\,r + 225\,p^{2} - 120\,m^{3}\,r \\
\mbox{} + 40\,m^{5}\,p + 265\,m^{2}\,p^{2} - 40\,q\,m^{3} - 80\,m
^{4}\,p - 20\,q\,m\,n + 360\,m^{2}\,p\,n + 30\,n^{2}\,m^{3} \\
\mbox{} + 600\,p\,q - 510\,m\,p^{2} - 120\,n^{3}\,m - 680\,m\,p\,q + 260\,m^{2}\,n\,q + 300\,m\,n\,r - 500\,n\,r \\
\mbox{} + 80\,p\,n^{2} - 170\,m^{3}\,p\,n + 60\,n^{3}))/(2\,m^{2} - 5\,n) }
}
\end{maplelatex}

\end{maplegroup}
\begin{maplegroup}
so alpha is just a number calculated directly from the coefficients of
the Quintic. Now the coefficient multiplying the linear term in d is
also linear in both eta and xi, solving it for eta and substituting
this into the zeroth term in d, gives a Quadratic Equation only in xi,
 i.e. let

\end{maplegroup}
\begin{maplegroup}
\mapleresult
\begin{maplelatex}
\mapleinline{inert}{2d}{zeroth_term_in_d := xi2*xi^2+xi1*xi+xi0;}{%
\[
\mathit{zeroth\_term\_in\_d} := \xi 2\,\xi ^{2} + \xi 1\,\xi  + 
\xi 0
\]
}
\end{maplelatex}

\end{maplegroup}
\begin{maplegroup}
with the equations for xi2, xi1, and xi0 given in the appendix(see the example given).  xi is
given by the Quadratic Formula,

\end{maplegroup}
\begin{maplegroup}
\mapleresult
\begin{maplelatex}
\mapleinline{inert}{2d}{xi := 1/2*(-xi1+sqrt(xi1^2-4*xi2*xi0))/xi2;}{%
\[
\xi  := {\displaystyle \frac {1}{2}} \,{\displaystyle \frac { - 
\xi 1 + \sqrt{\xi 1^{2} - 4\,\xi 2\,\xi 0}}{\xi 2}} 
\]
}
\end{maplelatex}

\end{maplegroup}
\begin{maplegroup}
Poly2 is cubic in d, setting it to zero, and using Cardano's Rule(on
Maple) to solve for d,  gives,

\end{maplegroup}
\begin{maplegroup}
\mapleresult
\begin{maplelatex}
\mapleinline{inert}{2d}{d :=
1/6*(36*d1*d2*d3-108*d0*d3^2-8*d2^3+12*sqrt(3)*sqrt(4*d1^3*d3-d1^2*d2^
2-18*d1*d2*d3*d0+27*d0^2*d3^2+4*d0*d2^3)*d3)^(1/3)/d3-2/3*(3*d1*d3-d2^
2)/(d3*(36*d1*d2*d3-108*d0*d3^2-8*d2^3+12*sqrt(3)*sqrt(4*d1^3*d3-d1^2*
d2^2-18*d1*d2*d3*d0+27*d0^2*d3^2+4*d0*d2^3)*d3)^(1/3))-1/3*d2/d3;}{%
\maplemultiline{
d := {\displaystyle \frac {1}{6}} (36\,\mathit{d1}\,\mathit{d2}\,
\mathit{d3} - 108\,\mathit{d0}\,\mathit{d3}^{2} - 8\,\mathit{d2}^{3} \\
\mbox{} + 12\,\sqrt{3}\,\sqrt{4\,\mathit{d1}^{3}\,\mathit{d3} - 
\mathit{d1}^{2}\,\mathit{d2}^{2} - 18\,\mathit{d1}\,\mathit{d2}\,
\mathit{d3}\,\mathit{d0} + 27\,\mathit{d0}^{2}\,\mathit{d3}^{2} + 4\,\mathit{d0}\,\mathit{d2}^{3}}\,\mathit{d3})^{(1/3)}/\mathit{d3} \\
\mbox{} - {\displaystyle \frac {2}{3}} (3\,\mathit{d1}\,\mathit{
d3} - \mathit{d2}^{2})/(\mathit{d3}(36\,\mathit{d1}\,\mathit{d2}
\,\mathit{d3} - 108\,\mathit{d0}\,\mathit{d3}^{2} - 8\,\mathit{d2
}^{3} \\
\mbox{} + 12\,\sqrt{3}\,\sqrt{4\,\mathit{d1}^{3}\,\mathit{d3} - \mathit{d1}^{2}\,\mathit{d2}^{2} - 18\,\mathit{d1}\,\mathit{d2}\,\mathit{d3}\,\mathit{d0}
 + 27\,\mathit{d0}^{2}\,\mathit{d3}^{2} + 4\,\mathit{d0}\,
\mathit{d2}^{3}}\,\mathit{d3})^{(1/3)}) \\
\mbox{} - {\displaystyle \frac {1}{3}} \,{\displaystyle \frac {
\mathit{d2}}{\mathit{d3}}}  }
}
\end{maplelatex}

\end{maplegroup}
\begin{maplegroup}
where d0, d1, d2, d3, and d4 are given in the example in the appendix.
With these substitutions the transformed Quintic, poly, takes the form
y\symbol{94}5+Ay+B =0 where A and B are generated by maple's determinant command(these are also in
the example given in the appendix).

\end{maplegroup}
\begin{maplegroup}
Then with a linear transformation the transformed equation
y\symbol{94}5+Ay+B = 0 , becomes z\symbol{94}5-z-s = 0, where

\end{maplegroup}
\begin{maplegroup}
\mapleresult
\begin{maplelatex}
\mapleinline{inert}{2d}{y := (-A)^(1/4)*z;}{%
\[
y := ( - A)^{(1/4)}\,z
\]
}
\end{maplelatex}

\begin{maplelatex}
\mapleinline{inert}{2d}{s := -B/((-A)^(5/4));}{%
\[
s :=  - {\displaystyle \frac {B}{( - A)^{(5/4)}}} 
\]
}
\end{maplelatex}

\end{maplegroup}
\begin{maplegroup}
We have used the fact that Poly4 is linear in a to get Poly4 = 0. Then
b,c and d were considered a point in space, on the curve of
intersection of a quadratic surface, Poly3, and a cubic
surface, Poly2. Giving Poly3 = Poly2 = 0, as required
[Cayley, 13],[Green, 14] and [Bring, 3].

\end{maplegroup}
\begin{maplegroup}
\begin{Heading 2}
THE SOLUTION TO BRING'S NORMAL FORM
\end{Heading 2}

\end{maplegroup}
\begin{maplegroup}
Bring's Normal Form is solvable. Any polynomial that can be transformed to
z\symbol{94}n-az\symbol{94}m-b=0
can be solved with the Hypergeometric Equation[Weisstein,5]. The solution is given by considering
z=z(s) and differentiating z\symbol{94}5-z-s = 0 w.r.t. s four times.
Then equate the fourth, third, second, first and zeroth order
differentials, multiplied by free parameters, to zero[Cockle, 8 and
9],[Harley, 10],[Cayley, 11]. Then make the substitution s\symbol{94}4 =
t. The resulting equation is a Generalized Hypergeometric Equation of
the Fuchsian type[Slater, 15], with the following solution
[Weisstein, 5],
\mapleinline{active}{1d}{                       }{%
}

\end{maplegroup}
\begin{maplegroup}
\mapleresult
\begin{maplelatex}
\mapleinline{inert}{2d}{z := -s*hypergeom([3/5, 2/5, 1/5, 4/5],[5/4, 3/4,
1/2],3125/256*s^4);}{%
\[
z :=  - s\,\mathrm{hypergeom}([{\displaystyle \frac {3}{5}} , \,
{\displaystyle \frac {2}{5}} , \,{\displaystyle \frac {1}{5}} , 
\,{\displaystyle \frac {4}{5}} ], \,[{\displaystyle \frac {5}{4}
} , \,{\displaystyle \frac {3}{4}} , \,{\displaystyle \frac {1}{2
}} ], \,{\displaystyle \frac {3125}{256}} \,s^{4})
\]
}
\end{maplelatex}

\end{maplegroup}
\begin{maplegroup}

Now calculating y, with

\end{maplegroup}
\begin{maplegroup}
\mapleresult
\begin{maplelatex}
\mapleinline{inert}{2d}{y := (-A)^(1/4)*z;}{%
\[
y := ( - A)^{(1/4)}\,z
\]
}
\end{maplelatex}

\end{maplegroup}
\begin{maplegroup}
We now undo the Tschirnhausian Transformation by substituting d, c, b, a
and y into the quartic substitution, Tsh1. The resulting Quartic
Equation is then solved using Ferrari's method[King, 6], this gives the
following set of equations,

\end{maplegroup}
\begin{maplegroup}
\mapleresult
\begin{maplelatex}
\mapleinline{inert}{2d}{g :=
1/12*(-36*c*d*b-288*y*c-288*a*c+108*b^2+108*a*d^2+108*y*d^2+8*c^3+12*s
qrt(18*d^2*b^2*y+18*d^2*b^2*a-3*d^2*b^2*c^2+576*d*b*a^2+576*d*b*y^2+76
8*y*a*c^2-432*y^2*c*d^2-432*y*c*b^2+1152*d*b*y*a+240*d*b*y*c^2+240*d*b
*a*c^2-54*c*d^3*b*a-54*c*d^3*b*y-864*y*c*a*d^2-432*a*c*b^2-2304*y^2*a+
12*y*d^2*c^3+12*d^3*b^3+12*a*d^2*c^3+162*a*d^4*y-432*a^2*c*d^2-48*a*c^
4+384*a^2*c^2-48*y*c^4-2304*y*a^2+384*y^2*c^2+81*y^2*d^4+81*a^2*d^4+12
*b^2*c^3+81*b^4-768*y^3-54*c*d*b^3-768*a^3))^(1/3)-12*(1/12*d*b-1/3*y-
1/3*a-1/36*c^2)/((-36*c*d*b-288*y*c-288*a*c+108*b^2+108*a*d^2+108*y*d^
2+8*c^3+12*sqrt(18*d^2*b^2*y+18*d^2*b^2*a-3*d^2*b^2*c^2+576*d*b*a^2+57
6*d*b*y^2+768*y*a*c^2-432*y^2*c*d^2-432*y*c*b^2+1152*d*b*y*a+240*d*b*y
*c^2+240*d*b*a*c^2-54*c*d^3*b*a-54*c*d^3*b*y-864*y*c*a*d^2-432*a*c*b^2
-2304*y^2*a+12*y*d^2*c^3+12*d^3*b^3+12*a*d^2*c^3+162*a*d^4*y-432*a^2*c
*d^2-48*a*c^4+384*a^2*c^2-48*y*c^4-2304*y*a^2+384*y^2*c^2+81*y^2*d^4+8
1*a^2*d^4+12*b^2*c^3+81*b^4-768*y^3-54*c*d*b^3-768*a^3))^(1/3))+1/6*c;
}{%
\maplemultiline{
g := {\displaystyle \frac {1}{12}} ( - 36\,c\,d\,b - 288\,y\,c - 
288\,a\,c + 108\,b^{2} + 108\,a\,d^{2} + 108\,y\,d^{2} + 8\,c^{3}
 + 12\mathrm{sqrt}( \\
18\,d^{2}\,b^{2}\,y + 18\,d^{2}\,b^{2}\,a - 3\,d^{2}\,b^{2}\,c^{2
} + 576\,d\,b\,a^{2} + 576\,d\,b\,y^{2} + 768\,y\,a\,c^{2} \\
\mbox{} - 432\,y^{2}\,c\,d^{2} - 432\,y\,c\,b^{2} + 1152\,d\,b\,y\,a + 240\,d\,b\,y\,c^{2} + 240\,d\,b\,a\,c^{2} \\
\mbox{} - 54\,c\,d^{3}\,b\,a - 54\,c\,d^{3}\,b\,y - 864\,y\,c\,a\,d^{2} - 432\,a\,c\,b^{2} - 2304\,y^{2}\,a + 12\,y\,d^{2}\,c^{3} \\
\mbox{} + 12\,d^{3}\,b^{3} + 12\,a\,d^{2}\,c^{3} + 162\,a\,d^{4}\,y - 432\,a^{2}\,c\,d^{2} - 48\,a\,c^{4} + 384\,a^{2}\,c^{2} - 48\,y\,c^{4} \\
\mbox{} - 2304\,y\,a^{2} + 384\,y^{2}\,c^{2} + 81\,y^{2}\,d^{4} + 81\,a^{2}\,d^{4} + 12\,b^{2}\,c^{3} + 81\,b^{4} - 768\,y^{3} \\
\mbox{} - 54\,c\,d\,b^{3} - 768\,a^{3}))^{(1/3)}\mbox{} - 12(
{\displaystyle \frac {1}{12}} \,d\,b - {\displaystyle \frac {1}{3
}} \,y - {\displaystyle \frac {1}{3}} \,a - {\displaystyle 
\frac {1}{36}} \,c^{2}) \left/ {\vrule 
height0.56em width0em depth0.56em} \right. \!  \! ( - 36\,c\,d\,b
 \\
\mbox{} - 288\,y\,c - 288\,a\,c + 108\,b^{2} + 108\,a\,d^{2} + 108\,y\,d^{2} + 8\,c^{3} + 12\mathrm{sqrt}(18\,d^{2}\,b^{2}\,y \\
\mbox{} + 18\,d^{2}\,b^{2}\,a - 3\,d^{2}\,b^{2}\,c^{2} + 576\,d\,b\,a^{2} + 576\,d\,b\,y^{2} + 768\,y\,a\,c^{2} - 432\,y^{2}\,c\,d^{2} \\
\mbox{} - 432\,y\,c\,b^{2} + 1152\,d\,b\,y\,a + 240\,d\,b\,y\,c^{2} + 240\,d\,b\,a\,c^{2} - 54\,c\,d^{3}\,b\,a \\
\mbox{} - 54\,c\,d^{3}\,b\,y - 864\,y\,c\,a\,d^{2} - 432\,a\,c\,b^{2} - 2304\,y^{2}\,a + 12\,y\,d^{2}\,c^{3} + 12\,d^{3}\,b^{3} \\
\mbox{} + 12\,a\,d^{2}\,c^{3} + 162\,a\,d^{4}\,y - 432\,a^{2}\,c\,d^{2} - 48\,a\,c^{4} + 384\,a^{2}\,c^{2} - 48\,y\,c^{4} - 2304\,y\,a^{2} \\
\mbox{} + 384\,y^{2}\,c^{2} + 81\,y^{2}\,d^{4} + 81\,a^{2}\,d^{4}
 + 12\,b^{2}\,c^{3} + 81\,b^{4} - 768\,y^{3} - 54\,c\,d\,b^{3} - 
768\,a^{3})) \\
^{(1/3)}\mbox{} + {\displaystyle \frac {1}{6}} \,c }
}
\end{maplelatex}

\end{maplegroup}
\begin{maplegroup}
\mapleresult
\begin{maplelatex}
\mapleinline{inert}{2d}{e := sqrt(1/4*d^2+2*g-c);}{%
\[
e := \sqrt{{\displaystyle \frac {1}{4}} \,d^{2} + 2\,g - c}
\]
}
\end{maplelatex}

\end{maplegroup}
\begin{maplegroup}
\mapleresult
\begin{maplelatex}
\mapleinline{inert}{2d}{f := 1/2*(d*g-b)/e;}{%
\[
f := {\displaystyle \frac {1}{2}} \,{\displaystyle \frac {d\,g - 
b}{e}} 
\]
}
\end{maplelatex}

\end{maplegroup}
\begin{maplegroup}
This gives four roots y1, y2, y3, and y4. These are then substituted
back into the Tschirnhausian Quartic to see which one satisfies it. The
root that satisfies it, will  also satisfy the General Quintic Equation,
Eq1[Prasolov and Solovyev, 8]. Let this root be r1. It is the root that
satisfies both the Quartic Tshirnhausian Transformation and General
Quintic Equation. This root varies as function of the parameters m, n, p,
q, and r. Another way of keeping track of which root satisfies Quintic and
Quartic, one could make a five column spread sheet or a five
dimensional plot of m, n, p, q and r. Once the root that satisfies the
quintic is determined, the other four roots of the quintic are
obtained by factoring out the root just obtained, this gives the
following equations,

\end{maplegroup}
\begin{maplegroup}
\mapleresult
\begin{maplelatex}
\mapleinline{inert}{2d}{r1 := yN;}{%
\[
\mathit{r1} := \mathit{yN}
\]
}
\end{maplelatex}

\end{maplegroup}
\begin{maplegroup}
\mapleresult
\begin{maplelatex}
\mapleinline{inert}{2d}{N*epsilon*[1, 2, 3, 4];}{%
\[
N\,\varepsilon \,[1, \,2, \,3, \,4]
\]
}
\end{maplelatex}

\end{maplegroup}
\begin{maplegroup}
\mapleresult
\begin{maplelatex}
\mapleinline{inert}{2d}{dd := m+r1;}{%
\[
\mathit{dd} := m + \mathit{r1}
\]
}
\end{maplelatex}

\end{maplegroup}
\begin{maplegroup}
\mapleresult
\begin{maplelatex}
\mapleinline{inert}{2d}{cc := n+r1^2+m*r1;}{%
\[
\mathit{cc} := n + \mathit{r1}^{2} + m\,\mathit{r1}
\]
}
\end{maplelatex}

\end{maplegroup}
\begin{maplegroup}
\mapleresult
\begin{maplelatex}
\mapleinline{inert}{2d}{bb := p+r1*n+r1^3+m*r1^2;}{%
\[
\mathit{bb} := p + \mathit{r1}\,n + \mathit{r1}^{3} + m\,\mathit{
r1}^{2}
\]
}
\end{maplelatex}

\end{maplegroup}
\begin{maplegroup}
\mapleresult
\begin{maplelatex}
\mapleinline{inert}{2d}{aa := q+r1*p+r1^2*n+r1^4+m*r1^3;}{%
\[
\mathit{aa} := q + \mathit{r1}\,p + \mathit{r1}^{2}\,n + \mathit{
r1}^{4} + m\,\mathit{r1}^{3}
\]
}
\end{maplelatex}

\end{maplegroup}
\begin{maplegroup}
Where x\symbol{94}4+ddx\symbol{94}3+ccx\symbol{94}2+bbx+aa=0 was
solved using Ferrari's method. This gives the other four roots of the
General Quintic Equation, r2, r3, r4 and r5. The last step to do is to
just check and make sure that they satisfy the original Quintic, and of
course they do! Now we have the five roots of the most general fifth
degree polynomial in a closed form. To write them all down on a piece
of paper, one would need, a piece of paper the size of large asteroid. But
these days with computers with such large memories, this can be done quite
easily. Now one might say what is the purpose of doing this? Why make
monstrous equations like this? Surely this must be useless? The
answer to these questions is quite simple! Was the quadratic equation
important and all the associated geometry of parabolas, circles
hyperbolas, ..etc? Were the cubic and quartic solutions important in
calculating arcs of ellipses, lemniscate curves, pendulum
motion, precession of planetary orbits, etc.. Now having the
equation for the roots of the quintic, we can investigate it's properties
and begin to use it to solve  physical problems. I think it is quite 
exciting that with the help of computer
algebra we can attack the non-linear problems of physics and
mathematics much more easily than in the days of Jacobi, Bring, Abel,
Cayley,..etc.   I hope that actually calculating the roots has dispelled
the common believe of most people I have talked to, that "it is impossible
to calculate the roots of the General Quintic Equation in a closed form". Now
I will put an end to this little project, by showing a maple session as
an example in the appendix. The other cases in the appendix(m=0 and n=0)
work as well, I do not want to waste time and space here by doing them..
 In the appendix I  calculate the roots of an arbitrary General
Quintic. I will let the reader load up the other equations on his/her
computer, and see for themselves that they, in fact, do work.   

\end{maplegroup}
\begin{maplegroup}
    These roots actually satisfy the Quintic identically, but the
computer ran out memory space. This equation which calculates the
roots of the Quintic, was checked with various values of the
parameters m, n, p, q and r. It works for all values except one. When m =
0. This is probably because all the equations above really need to be
put into one another then the division by zero cancels. So below, I
diveded the calculation into three cases only to avoid having the
computer divide by zero. When all the equations are put together, I get a
Maple error saying "Object Big Prod". Maple probably has a memory
protection limit built into it. Once this is removed, then run on a
computer with a larger memory, then all the above equations may be
substituted into one another. Also the final equation with all the
substitutions completed, may be decrease in size, due to cancellations.
The creators of Maple tell me that this is going to be accomplished with
the next version of Maple, Maple VI.     

\end{maplegroup}
\begin{maplegroup}
\begin{mapleinput}
\end{mapleinput}

\end{maplegroup}
\begin{maplegroup}
\begin{Heading 2}
{\large  }APPENDIX
\end{Heading 2}

\end{maplegroup}
\begin{maplegroup}
\begin{mapleinput}
\mapleinline{active}{1d}{             EXAMPLE: A MAPLE SESSION}{%
}
\end{mapleinput}

\end{maplegroup}
\begin{maplegroup}
\begin{mapleinput}
\mapleinline{active}{1d}{restart;}{%
}
\end{mapleinput}

\end{maplegroup}
\begin{maplegroup}
\begin{mapleinput}
\mapleinline{active}{1d}{Digits := 200:}{%
}
\end{mapleinput}

\end{maplegroup}
\begin{maplegroup}
\begin{mapleinput}
\mapleinline{active}{1d}{m := -200*I;}{%
}
\end{mapleinput}

\mapleresult
\begin{maplelatex}
\mapleinline{inert}{2d}{m := -200*I;}{%
\[
m :=  - 200\,I
\]
}
\end{maplelatex}

\end{maplegroup}
\begin{maplegroup}
\begin{mapleinput}
\mapleinline{active}{1d}{n := 1340;}{%
}
\end{mapleinput}

\mapleresult
\begin{maplelatex}
\mapleinline{inert}{2d}{n := 1340;}{%
\[
n := 1340
\]
}
\end{maplelatex}

\end{maplegroup}
\begin{maplegroup}
\begin{mapleinput}
\mapleinline{active}{1d}{p := 1.23491*10^1;}{%
}
\end{mapleinput}

\mapleresult
\begin{maplelatex}
\mapleinline{inert}{2d}{p := 12.34910;}{%
\[
p := 12.34910
\]
}
\end{maplelatex}

\end{maplegroup}
\begin{maplegroup}
\begin{mapleinput}
\mapleinline{active}{1d}{q := -2.39182*10^2;}{%
}
\end{mapleinput}

\mapleresult
\begin{maplelatex}
\mapleinline{inert}{2d}{q := -239.18200;}{%
\[
q := -239.18200
\]
}
\end{maplelatex}

\end{maplegroup}
\begin{maplegroup}
\begin{mapleinput}
\mapleinline{active}{1d}{r := 3.3921817*10^2;}{%
}
\end{mapleinput}

\mapleresult
\begin{maplelatex}
\mapleinline{inert}{2d}{r := 339.2181700;}{%
\[
r := 339.2181700
\]
}
\end{maplelatex}

\end{maplegroup}
\begin{maplegroup}
\begin{mapleinput}
\mapleinline{active}{1d}{To avoid having the computer divide by zero, the calculation of
alpha, eta, and xi is divided into three cases: 1) m and n not equal to zero 2) m equal to zero and n not 3) both m and n equal to zero. The
last one is Bring's original transformation.}{%
}
\end{mapleinput}

\end{maplegroup}
\begin{maplegroup}
\begin{mapleinput}
\mapleinline{active}{1d}{m and n not equal to zero}{%
}
\end{mapleinput}

\end{maplegroup}
\begin{maplegroup}
\begin{mapleinput}
\mapleinline{active}{1d}{alpha :=
evalf(1/2*(-13*m*n-10*n^2+4*m^3+20*q+17*m^2*n-4*m^4+15*p-17*m*p+sqrt(-
680*q*m*p+200*m^2*r+30*m^3*n^2+360*p*m^2*n-80*m^4*p-500*n*r+260*q*m^2*
n-190*m*n*p-80*m*p*n^2-15*n^2*m^4+60*n^3*m^2+80*p*n^2+60*n^3-170*m^3*n
*p-200*n*q+80*m^2*q+40*m^5*p-100*q*n^2+400*q^2-120*m^3*r+225*p^2+265*m
^2*p^2+300*m*n*r-40*q*m^4+60*n*p^2+600*p*q-15*m^2*n^2+40*m^3*p-510*m*p
^2-120*m*n^3-40*m^3*q-20*m*n*q))/(2*m^2-5*n)):}{%
}
\end{mapleinput}

\end{maplegroup}
\begin{maplegroup}
\begin{mapleinput}
\mapleinline{active}{1d}{xi3 :=
evalf(580*m^3*p*alpha^2*n-80*m^5*p*alpha^2+1500*m^3*p^2*alpha+5200*m^2
*n^3*q+1360*m^6*n*q+1600*q^2*m^2*alpha+5000*q^2*m*n-4000*q^3+3200*q*n^
3*alpha-320*q*m^6*alpha-5775*m*p^2*n*alpha-4065*m^4*n^2*q-6625*m^2*n*q
^2-5625*q*m*p^2+3285*m^2*p^2*q+5820*m^4*p^2*n-160*q*m^4*alpha^2-8020*m
}{%
}
\end{mapleinput}

\end{maplegroup}
\begin{maplegroup}
\begin{mapleinput}
\mapleinline{active}{1d}{
^2*p^2*n^2-1580*m^4*p^2*alpha+310*m*p*n^4+3300*m*p*q^2+860*m^7*p*n-360
*m^5*p*q+5895*q*m^3*n^2-4000*q^2*n*alpha+1040*q*m^4*n-1990*q*m^2*n^2-2
400*q*m^5*n-180*n^5+200*m*p*q*n^2+1125*n*p^2*alpha^2-2250*n^2*p^2*alph
a-375*n^2*m^2*r-1500*n^3*m*r+1000*n^2*p*r+375*n^2*m^3*r-5000*n*q*r+495
0*m^2*p^3-400*q*m*p*n*alpha+760*q*m^4*p+320*q*m^5*alpha-2820*m^5*p*n^2
+2585*m^3*p*n^3+3800*q*n^2*p-2000*q^2*m^3-5300*q*m*n^3+2250*q*p^2-2250
*m*p^3+30*n^2*m^4*alpha^2-160*q*m^6-850*m^4*p^2+160*m^8*p-80*m^7*p+104
5*n^3*p^2-2000*n*q^2-3125*n*r^2+1500*n^3*r+1200*n^3*q-400*q*m^3*p+60*n
^2*m^6*alpha+7055*m^2*p^2*n*alpha-1485*n^4*m^3+525*n^4*m^2+800*m^2*q^2
+1250*m^2*r^2-240*n^3*m^4+540*n^3*m^5+5625*p^2*r+1780*m^5*p^2-195*n^3*
m^2*alpha^2-675*n^2*p^2+300*n^4*alpha^2+320*q*m^7-600*n^5*alpha-2005*n
^3*m^2*p+800*q*m^2*alpha^2*n+435*n^3*m^3*alpha-780*n^4*m*alpha+2600*q*
m*n^2*alpha+2000*m^2*q*r+1000*m^3*r*p-450*m^2*p^2*alpha^2-4275*p^3*n-6
0*n^2*m^7-1140*n^4*p-1000*m^4*p*r+3375*p^3*alpha-950*m*p*n^2*alpha^2+1
140*n^5*m-500*m^3*q*r+30*n^2*m^6+7885*n^2*m*p^2-160*m^8*q+5200*n^2*q^2
-60*n^2*m^5*alpha-2200*n^4*q+4230*n^2*m^4*p+1650*m^4*q^2+30*m^8*n^2-11
55*m^2*n^5+7500*q*p*r+990*m^4*n^4-300*m^6*n^3+4500*q*p^2*alpha-5700*q*
p^2*n-1840*q*m^3*n*alpha+280*m^3*p*q*n+900*n^3*alpha*p-6375*m*p^2*r-38
25*m*p^3*alpha+4845*m*p^3*n-1170*q*m^2*n*p-30*n^3*m*p-3000*q*n*alpha*p
-560*q*m^3*p*alpha-160*m^7*p*alpha-495*n^3*m^4*alpha+1170*n^4*m^2*alph
a-1490*n^2*m^3*p+100*q*m*n*p+160*m^6*alpha*p+700*m^5*n*p-4750*n*r*m*p-
1560*m^6*n*p+4500*n*m^2*p*r-250*n*m*q*r+1200*q*m^2*alpha*p-3375*n^2*m^
3*p*alpha+3700*n*m^2*p^2+2160*q*m^4*n*alpha-4500*q*m^2*n^2*alpha-930*m
^6*p^2-9440*n*m^3*p^2+880*n^3*m*p*alpha+2245*n^2*m^2*alpha*p-80*m^9*p+
1440*m^5*p*n*alpha-2720*m^3*p^3-1280*m^4*n*alpha*p-1000*q*n^2*alpha^2+
300*n^6):}{%
}
\end{mapleinput}

\end{maplegroup}
\begin{maplegroup}
\begin{mapleinput}
\mapleinline{active}{1d}{xi2 :=
evalf(5640*m^6*alpha*q*p+46000*m*q^2*r-240*m^6*alpha*p^2+3400*m^4*alph
a*q*r-480*m^8*alpha*q+4875*p^2*alpha*n*r-990*m^5*alpha*n^4+1800*m^4*al
pha^2*q*p-600*n^4*alpha^2*p+4750*n^2*alpha^2*p^2*m+160*m^8*alpha*r-191
0*n^3*alpha^2*p*m^2+240*m^7*alpha^2*q-4860*m^5*alpha*q^2+120*m^5*alpha
^2*p^2+60*m^5*alpha^2*n^3-15000*p*alpha*q*r-4680*m^7*alpha*q*n-7600*p^
2*alpha*m^2*r-3900*n^3*m*alpha*r+240*m^7*alpha*p^2+1500*n^3*alpha^2*r+
12600*q^2*alpha*m*p-160*n*alpha^2*r*m^4+600*n^5*alpha^2*m-19250*m*n^2*
r^2+27500*n*p*q*r+4500*p*alpha*n^2*r+5100*p^2*alpha*q*n-200*m^6*alpha^
2*n*p-2220*n^3*p*q-8700*n^3*p*r+3420*m^6*alpha*n^2*p+26250*n*p*r^2+340
0*n*p*q^2-14750*m^2*p*r^2-400*m^8*alpha*n*p-10980*m^4*alpha*q*n^2+280*
m^7*p*r-2250*p^3*alpha^2*n-520*m^5*alpha*r*p-1800*p^2*alpha*n^3+900*p^
3*alpha^2*m^2-9375*p*alpha*r^2+14820*m^5*alpha*q*n^2+6000*r*q*n^2+1600
*r*q*m^4-6520*m^5*alpha*n*p^2+2460*m^2*n^3*r-5400*m^5*alpha*q*p+400*m^
7*alpha*n*p-500*m^2*n*r^2-3020*m^5*alpha*p*n^2-280*m^6*n*r-640*m^4*n^2
*r-4560*m*n^5*p+8700*m*n^4*r-120*m^6*alpha*n^3-600*m^6*alpha*n*r-1920*
m^5*alpha^2*q*n+1360*m^4*alpha^2*n^2*p-840*m^4*alpha*n^2*r-7320*m^4*al
pha*n^3*p+870*m^4*alpha*n^4+3360*m^4*alpha*q^2-3720*m^7*p*n^2+4500*n^2
*alpha^2*m^3*q+14050*m*n^3*p*r+1700*n^2*alpha^2*q*p-10095*p^2*m*q^2-10
920*q^2*alpha*m^2*n-700*n^2*alpha^2*m^2*r+11000*m*p*q*r-2200*n*alpha^2
*m^3*p^2-3380*m^5*p^2*n-4350*p^3*m^2*q-870*m^2*p*n^4-2560*m^3*p*n^4-30
00*n^3*alpha^2*m*q-4000*m^3*alpha*q*r+4200*n*alpha^2*m*q^2-1000*m^2*p^
2*r-5000*n*alpha^2*q*r+7120*m^2*p*q*n^2+6080*m^4*p^2*n*alpha+830*m^2*p
^2*q*n-18300*m^4*p^2*q-47000*m^2*p*q*r-38890*m^3*p*q*n^2-13100*m^2*p^3
*n+12080*m^3*p^2*n^2+23580*m^2*p^2*n^3+31040*m^3*p^3*n+33660*m^5*p*q*n
}{%
}
\end{mapleinput}

\end{maplegroup}
\begin{maplegroup}
\begin{mapleinput}
\mapleinline{active}{1d}{
-13480*m^4*p*q*n+21850*m*p^2*n*r+820*m^4*p*n*r-8200*m^2*n*q*r-3300*m*p
*n^2*r+200*m^5*p*r+4660*m^3*p*n^3*alpha-480*m^6*p*r-18800*m^2*p^2*n^2*
alpha-13500*m^2*p^2*q*alpha-200*m^4*p*r*alpha+6100*m*p^2*q*n+1500*m*p^
2*alpha*r+1000*m^3*p*n*r-1870*m^2*p*n^2*r-2700*m^3*p^3*alpha-17250*q*a
lpha*m^3*n^3+2500*m*p*r^2+4500*m*p^4-25855*q^2*alpha*m^2*p+22705*q^2*a
lpha*m^3*n+20100*m*p^3*n*alpha-4500*q*p^3+700*m^3*p^2*r+41940*m^2*q^2*
n*p-20200*m^4*q*n*r-31500*m^4*p^2*n^2+8400*m^3*p^2*q+7200*m^6*p^2*n-17
400*m^2*p*q^2+40980*m^3*p*q^2+9200*q^3*p-1650*m*p^2*n^3-26600*m*p^3*n^
2+24275*m^2*q*n^2*r-18745*m^2*q*n^3*p+9500*q*alpha*n^2*r-10000*q^2*alp
ha*r+14300*q*alpha*m^3*p^2-5900*q*alpha*n^3*p+7800*q*alpha*m*n^4-1680*
q^2*alpha^2*m^3-16300*q^2*alpha*m*n^2+7800*q*alpha*n^3*m^2+8400*q^3*al
pha*m-8860*n^2*q*p^2*m+8550*p^4*n-6750*p^4*alpha-11250*p^3*r+21780*m^5
*q^2*n-27875*m^2*q^2*r-2760*m^9*q*n-24630*m^4*q^2*p+3200*m^6*q*r-20420
*m^6*q*n*p+33390*m^4*q*n^2*p+10260*m^5*q*p^2-9010*m^3*q*n*p^2+3240*m^8
*q*p+25700*m^3*q*r*p+8000*q^2*alpha*n*p-17880*m^6*q*n^2+25290*m^4*q*n^
3+5040*m^8*q*n+30960*m^2*q^2*n^2+1350*n^2*p^3-34045*m^3*q^2*n^2+16350*
m^3*q*n^4-10000*r*q^2-120*m^6*n^2*r+1460*m^3*alpha*n^2*r-900*r*n^4-250
00*r^2*q-1000*r^2*m^4+2530*m^2*n^5*p+440*m^5*alpha*n*r+7500*r^2*n^2+23
40*m^3*n^5*alpha+2060*m^8*n^2*p+600*n^5*r-200*m^3*alpha^2*p*r-600*n^6*
p+37000*m^3*q*n*r+19740*m^5*n^2*p^2-41500*m*q*n^2*r-24900*m^3*n^3*p^2+
4730*m^4*n^4*p-6360*m^6*n^3*p+4200*n^4*q*p-2700*n^3*q*r-10600*n^2*q^2*
p+220*m^2*n^4*p*alpha+24010*m^3*n^2*p^2*alpha-12200*n*q^3*m+3000*n*q^2
*r-620*m^3*n^2*r*p-4800*n^5*q*m+13300*n^3*q^2*m+20860*n^3*q*m*p-24460*
n*q^2*m*p-33300*m^4*q^2*n-15625*r^3-14400*m^2*q*n^4-6460*p^2*q*n^2-160
*m^7*alpha*r+11250*p^3*q*m-16125*p^2*q*r-3820*p^2*m^7*n+3360*n^4*alpha
*m*p-360*m*n^6+1870*p^2*m*n^4-9675*p^3*q*alpha+12255*p^3*q*n+500*m^3*a
lpha*r^2+3875*n^3*alpha*m^2*r+80*m^6*alpha^2*r+600*n^7*m+6840*m^5*q*n^
2-7080*m^3*q*n^3+2000*m^2*alpha^2*q*r+4200*m^6*alpha*q*n-2280*m^7*q*n-
900*n*m^5*r*p-8520*m*q^2*n^2+13000*n*m*alpha*q*r+13320*m^3*q^2*n+2280*
m^2*n^6+3000*m*q*n^4+29375*m*q*r^2-6000*m^5*q*r-24700*p^3*m^2*n*alpha-
600*m^7*n^4-7260*p^2*n^3*m*alpha+4210*m^4*n^3*r+1980*m^5*n^5-7075*m^2*
n^4*r+60*m^9*n^3+80*r*m^8-17290*p^2*m^2*n*r+360*n^5*p+11280*m^7*q*n^2-
440*n*m^8*r-200*n*m^10*p+4500*p^3*n^2*alpha+3020*p^3*alpha*m^4+9000*p^
4*alpha*m+1020*p^2*r*m^4-8175*p^2*n^2*r+15000*p^3*r*m-11400*p^4*n*m+29
390*p^3*n^2*m^2-18360*p^3*m^4*n-13420*q*alpha*m*p*n^2+3000*m^6*p*q-624
0*m^7*p*q-200*m^8*p*n+400*m^9*p*n-3200*m^4*p*n^3+1660*m^6*p*n^2+9360*m
^5*p*n^3-160*m^9*r-390*n^4*alpha^2*m^3+80*m^10*r-1000*m^5*r^2-2880*m^5
*q^2-2090*p^3*n^3+480*m^9*alpha*q+5380*p^4*m^3+1620*p^3*m^6+120*m^7*al
pha*n^3+1500*m^4*p^3-26580*m^4*alpha*q*n*p-1560*n^5*alpha*m^2+240*m^11
*q-7130*m^3*n^3*r+1080*m^6*n^4+10500*m^3*n*r^2+8400*m*q^3+720*m^7*n*r+
840*m^5*n^2*r-5180*n*alpha^2*m^2*q*p-1200*n^6*alpha*m+1200*n^5*alpha*p
-2100*n^4*alpha*r+625*n*m*alpha*r^2+120*m^7*p^2-240*m^8*p^2+500*n*alph
a^2*r*m*p+6360*m^6*q^2+12525*m^3*q^3-1600*m^2*p*n*alpha*r-3120*m^5*p^3
-17250*q*alpha*m^2*n*r+22960*m^3*p*q*n*alpha+26750*q*alpha*r*m*p-2310*
m^3*n^6-120*m^8*n^3-3480*m^7*q^2-480*m^10*q-21000*m^2*q^3-20130*m^5*q*
n^3+6140*m^3*alpha*n*p*r+240*m^9*q-2970*m^4*n^5-11750*n*q*r*m*p+7820*q
*alpha*n*p^2*m+19505*q*alpha*n^2*p*m^2-9900*p^4*m^2+120*p^2*m^9+2280*p
^2*n^4+1050*m^3*n^5-7750*n^2*alpha*r*m*p+60*m^7*n^3-480*m^5*n^4):}{%
}
\end{mapleinput}

\end{maplegroup}
\begin{maplegroup}
\begin{mapleinput}
\mapleinline{active}{1d}{xi1 :=
evalf(6850*m^2*p^3*r*alpha+38040*q*p*m^5*n^3-5745*q*m^3*n^2*alpha^2*p-
14500*n*m^2*p^2*r*alpha+14550*m*p*n^2*r^2-13790*q^2*m^2*n^3*alpha+1754
0*q^2*m*n^2*p*alpha+4880*q*p*m^9*n+9000*q*m^2*r^2+2720*q*m^9*r-240*q*m
^11*p-2000*q*m^2*n^4*alpha^2+7160*q*m^8*n^3+5315*q^2*m^2*alpha^2*n^2+1
9815*q^2*m^4*n^2*alpha+13500*m^2*p*r^2*alpha-5625*q*p^4*m-12265*q*m^6*
n^4-2200*m^4*p^3*r-3040*q*m^8*n^2*alpha+320*q*m^10*n*alpha+24550*m^3*p
*r^2*n-3610*q^2*m^4*n*alpha^2-27750*m^3*p*r^2*alpha+160*q*m^12*n-500*q
*n^2*r^2+10710*q^3*m^4*n+8835*q^2*m^2*n^4-6250*q*r^2*alpha^2-14580*q*m
^5*p^3+25000*q*m^4*r^2+23720*q^2*m^2*n*p*alpha-43230*q*p^2*m^4*n^2+191
90*m^2*p^3*n*r-1240*q*m^6*n^2*alpha^2-3860*q*m^8*p^2-18600*m^5*p*r^2-5
000*m*p*r^2*alpha^2+2700*q*n^3*p*alpha^2*m+1500*m*p^3*r*alpha+160*q*m^
8*n*alpha^2-45030*q^2*p^2*m^2*n-40620*q^2*p*m^5*n+300*m^5*p^3*alpha^2-
300*m^3*p^3*r+6035*q^2*p^2*n^2+10050*q*p^2*m^2*n^3-13905*q*p*m^2*n^2*r
-7550*q*p*m^4*n*r-3980*q^2*m^6*n*alpha+7040*q*m^3*n^2*alpha*r+60*q*m^3
*n*alpha^2*r-20000*q*m^5*n*alpha*r+410*q^2*m^6*alpha^2+15340*q^2*m^6*n
^2+17645*m^2*p^4*n*alpha+9400*n*m^2*p^4+4950*m^2*p^5-2070*q*m^4*p^2*al
pha^2+9600*n*m^3*p^2*r+8360*q*n*p*alpha*m^7+20750*m*p*n*r^2*alpha-6680
*q*m^6*p^2*alpha+230*q^3*m^4*alpha+17310*q*p^3*m^3*n-21620*q*p*m^3*n^4
-3800*m*p^3*n^2*alpha^2-23980*q*p*m^7*n^2+2730*q*m^5*n*alpha^2*p+30940
*q*p^2*m^6*n-2250*m*p^5-1800*q*m^10*n^2-9000*m^3*p*r^2+570*m^3*p^3*alp
ha^2*n+80*q^2*m^7*alpha+34800*m^4*p*r^2+4320*q^2*m^7*n-14600*q*m*n*p*r
*alpha+28600*q^2*m^2*r*p+24660*q*m^4*n^2*alpha*p+4600*q*m^2*n^5*alpha+
31455*q^2*m^2*p^2*alpha+14160*q*p*m^2*n*r*alpha-700*q*n^2*p^2*alpha^2-
1800*q*m^2*p*alpha^2*r-80*q^2*m^8*alpha-11140*q^3*m^2*n^2-480*q*m^9*p*
alpha+20820*q*m^4*p*alpha*r+5170*q^2*m^5*p*alpha-240*q*m^7*p*alpha^2-4
050*q^2*m^3*n*r-11160*q^3*m*alpha*p+10500*q^2*m*alpha^2*r-15050*q^2*m^
3*alpha*r-14325*m*p^4*n*alpha-19700*q^2*m*alpha*n*r-7500*q*m*n^2*alpha
^2*r-11550*q^2*m^5*r-2000*q*p^2*n^4+3640*q^2*m^7*p-600*m^3*p^4*alpha-5
050*q^2*n*p^2*alpha+11500*q*m*n^3*alpha*r+15930*q^3*m*n*p-6745*q*n*p^3
*alpha*m+28895*q^2*m^4*p^2+6350*q^2*m*n^2*r-14570*q^2*m*n^3*p-7300*q*n
^4*p*alpha*m-23780*q^2*m^4*n^3-2320*q^2*m^8*n+19920*q*m^5*n^2*r-49380*
q^2*p^2*m^3-48590*q^2*p*m^3*n*alpha-5600*q^2*n^3*m^2+320*q^3*m^3*alpha
-23800*q^2*m^5*n^2-21020*q*m^3*r*p^2+14720*q*m^6*r*p-9750*q*m^2*r^2*n+
33250*q*m^2*r^2*alpha-20900*m*p^3*n*r+4250*q*n*p*alpha^2*r-8700*q*n^2*
p*alpha*r+8465*q*n*r*m*p^2+7245*q^2*p^3*m-16510*q^3*m^3*p-1000*m^2*p^3
*r-1400*q^2*p^2*n-380*q^2*m*n^4-28400*q*r*alpha*m*p^2-39100*m^2*p*r^2*
n+2700*q*n^3*p^2*alpha-5350*q^2*n*r*p+2855*q*m^4*n^3*alpha^2+3105*q*m^
2*p^2*alpha^2*n-3300*m^3*p^2*alpha^2*r-11200*q^3*m*p-47120*q^2*m^2*n^2
*p+2020*q^3*m*n^2-4350*q^2*m*alpha^2*n*p-10280*q*m^4*n^4*alpha-15760*q
*m^7*n*r+4640*q*m^7*r*alpha+21000*q^2*p^2*m^2-3100*q*n^4*r*m+8780*q*p^
3*n^2*m+2500*q*n^3*r*p+4600*q*p*n^5*m-675*q*m*p^3*alpha^2-22750*q*m*p*
r^2+1020*q*m^5*alpha^2*r+3750*q*n*r^2*alpha+600*m^4*p^4*alpha-7720*q*m
^6*n*p*alpha+5280*q^2*m^3*p*alpha^2+23260*q^2*m*n*p^2+9345*q*m^4*n^5-2
0220*q*m^3*p^3*alpha+13000*q^2*r*alpha*p+10220*q*m*n^2*p^2*alpha+14110
*q^3*m^2*alpha*n-4410*q^3*m^2*alpha^2-32940*q*m^3*n*alpha*p^2-6720*q^2
*m^6*p-18600*q^3*m^3*n+8760*q*m^4*n*alpha*r-12380*q*m^2*n^3*p*alpha+26
0*q^2*m*alpha*n^3-15600*q^2*m*alpha*p^2+40520*q*m*n^2*p*r+400*q*m^2*n*
p*r-585*q*m^3*n^3*r+8920*q^2*m*n^2*p-3860*q^2*m^4*p*alpha-16720*q*p*m^
3*n*r-27520*q^2*m^3*n*p-8920*q^2*m^3*n^2*alpha+2460*q^2*m^5*n*alpha+21
820*q^2*m^3*n^3-300*q^2*n^2*p*alpha-8200*q*m^3*p*alpha*r+24000*q^2*m^2
*n*r-1040*q^3*m*alpha*n-2800*q^2*m*r*n+66380*q^2*p*m^4*n+4280*q*m^6*n^
}{%
}
\end{mapleinput}

\end{maplegroup}
\begin{maplegroup}
\begin{mapleinput}
\mapleinline{active}{1d}{
3-30460*q*p*m^5*n^2*alpha+20100*q*p*m^2*n^4-9120*q*p*m^8*n-4220*q*p^3*
m^2*n+3280*q*m^9*n^2+380*q^2*n^3*p-2600*q*m^2*n^6+800*q^4-46000*q^2*m*
r*p+4420*q*m^3*n^4*alpha-320*q*m^9*n*alpha+3640*q*m^2*n^2*alpha*r+2660
*q*p^3*n^2-2000*q^2*m^6*n+1620*q*n^5*m^2+39570*q*p^2*m^3*n^2+2720*q*m^
7*n^2*alpha+15900*q*m^2*p^3*alpha+480*q*m^8*p*alpha-44220*q*m^4*n^3*p+
15155*q*n^4*m^5-3520*q*m^6*p^2-4130*q*n^4*m^4-2000*q^4*m+1250*q^2*r^2-
2500*q^3*m*r+1200*q^3*p*alpha+450*q^2*p^2*alpha^2+2000*q^3*r+1250*q^4*
m^2-6200*q*m*n*p^3-2880*q*m*n^4*p+4240*q*m^7*n*p-2200*q*m^2*n^2*p^2+10
45*n^3*p^4-6820*q*m^5*n^3*alpha-1125*p^3*r*alpha^2-6000*p*r^2*n^2+5370
*p^3*r*n^2+2500*p^2*r^2*m-2025*p^3*r*n*alpha-4800*p*r*q*n^2+12500*p*r^
3+8000*p*r*q^2+720*p*r*n^4+20000*p*r^2*q-18375*p^2*r^2*n+5375*p^2*r^2*
m^2+5625*p^2*r^2*alpha+38580*q*p*m^6*n^2-51320*q*p^2*m^5*n-14800*q*m*n
^3*p^2+6040*q*m^5*p^2*alpha+9200*q*m^3*n^3*p-15000*q*m^5*n^2*p-5040*q*
m^8*r+480*q*m^10*p+2320*q*m^7*r+160*q*m^10*n-240*q*m^9*p-1480*q*m^8*n^
2+5840*q*m^3*n^2*r+25820*q*p^3*m^4-7960*q*n^5*m^3-2100*q*n*p^3*alpha-3
20*q*m^11*n-11600*q*p^3*m^3+1125*n*p^4*alpha^2-2250*n^2*p^4*alpha-1682
0*q*n^3*m^2*r+2040*q*n^3*m*r-22345*n^2*m^2*p^4-18700*q*p^2*n*r+12000*q
*p^2*r*alpha+41600*q*p^2*m^2*r-7000*q*p^2*m*r-3040*q*m^6*r*alpha+24080
*q*m^6*n*r-12840*q*m^4*n^2*r+40460*q*p^2*m^4*n*alpha+22725*q*p*m^3*n^3
*alpha-10675*q*p^2*m^2*n^2*alpha+21160*q*m^4*n*p^2+7380*q*p^2*m^7-1128
0*q*n^3*m^7-195*n^5*m^4*alpha^2+9620*q^2*m^4*n^2+1020*q*p^2*n^3-10320*
q*m^5*r*n-7600*q^2*m^3*r-41500*q*m^3*r^2+3080*q^2*m^5*p+60365*q^2*p*m^
3*n^2+19500*q^2*m^4*r+7440*q^3*m^2*n+27760*q^3*m^2*p-500*q^2*n^2*r-152
0*q^3*n*p-240*n^2*m^9*p*alpha-1565*n^2*m^4*p^2*alpha^2+3830*n^3*m^2*p^
2*alpha^2-495*n^5*m^6*alpha-16695*n^4*p^2*m^4+12000*q*m^4*p*r-25520*q*
m^5*p*r-7000*q*m*r^2*n-440*q^3*n^2-1900*n^2*m^8*p^2+30*n^4*m^6*alpha^2
-960*n^4*m^3*alpha^2*p-120*n^2*m^11*p-8360*n^2*m^9*r+30*n^4*m^10+60*n^
4*m^8*alpha+300*n^6*m^2*alpha^2+1170*n^6*m^4*alpha+300*n^4*p^2*alpha^2
-2700*n^2*m^2*r^2+29350*n^3*p^3*m^3+2220*n^6*p*m^3-3540*n^4*p*m^7-600*
n^5*p*alpha^2*m+13160*n^3*p^2*m^6+1980*n^3*p*alpha*m^7-3900*n^2*m^6*p^
2*alpha+18715*n^2*p^4*m+70*n^5*p^2*m^2+8205*n^3*p*m^2*r*alpha+14180*n^
3*m^5*alpha*r+1875*n^2*r^2*alpha^2-600*n^5*p^2*alpha+1250*m*n*r^3+2145
*n^5*p*m^5-16740*n^2*m^5*p^3+1200*n^3*p*m^9+780*n^3*m^5*alpha^2*p+1576
0*n^3*m^7*r-120*n^2*m^7*p*alpha^2+1500*n^4*m*alpha^2*r-2100*n^5*m*alph
a*r+6380*n^3*p^3*alpha*m+4090*n^4*p*m^2*r-24850*n^2*m^4*r^2+30300*n^3*
p*m^4*r-3535*n^4*m^3*alpha*r+1835*n^3*m^3*alpha^2*r-24960*n^2*m^3*p^3*
alpha+345*q*m^2*p^4+9720*n*m^4*p^4+1200*n^6*p*alpha*m-6495*n^2*m^2*p*a
lpha^2*r-600*n^7*m^2*alpha-160*n^5*m^3*r-1500*n^3*p*alpha^2*r+4920*n^3
*m*p*r*alpha-9300*n^3*r*m*p^2-600*n^5*r*p+7200*n*m^4*r^2-2845*n^2*r*al
pha*m*p^2-11175*n^4*m^5*r+2100*n^4*p*alpha*r-21775*n^2*m^2*r^2*alpha-3
8820*n^2*m^6*r*p-2180*n^4*p^3*m-780*n^6*m^3*alpha-600*n^7*p*m-3410*n^2
*m^5*alpha^2*r+2415*n^4*m^4*alpha*p-11280*n^2*m^7*r*alpha-29890*n^2*m^
3*r*p^2+600*n^6*r*m-1500*n^3*r^2*alpha+2425*n^3*m^2*r^2+9300*q*m^6*n^3
*alpha+2460*n^5*m^2*p*alpha-2580*n^4*m*p^2*alpha-13370*n^3*m^3*alpha*p
^2+720*m^5*n*p^3-13750*m^3*r^3+2250*q*p^4-9840*n^4*m*p*r-3420*n^6*p*m^
2-7060*n^3*m^4*alpha*r-480*n^3*m^2*p*r-5290*q^3*p^2-40*q^2*m^10+19180*
n^2*m^3*p*alpha*r-1740*n^3*m^6*p*alpha-3945*n^4*p*m^5*alpha+15350*n^4*
p^2*m^3-2160*n^3*p*m^8+22500*m^2*r^3-1155*n^7*m^4-8625*p^4*m*r-300*n^5
*m^8-6555*p^4*q*n+5175*p^4*q*alpha+8625*p^3*q*r+6555*p^5*m*n-5175*p^5*
m*alpha+140*q^3*m^6+60*q^2*n^4-40*q^2*m^8+80*q^2*m^9+40*q^3*m^4-180*q^
3*m^5-26300*n^3*p*m^3*r+225*m^2*p^4*alpha^2-630*n^4*m^2*p^2-840*n^5*m^
}{%
}
\end{mapleinput}

\end{maplegroup}
\begin{maplegroup}
\begin{mapleinput}
\mapleinline{active}{1d}{
3*p-1710*n^4*m^5*p+3420*n^5*m*p^2-21580*n^3*p^2*m^5+5130*n^4*p*m^6+525
*n^6*m^4+900*n^3*p^3*alpha+240*n^2*m^8*p*alpha+21430*n^2*m^2*p^3*alpha
+300*n^6*p^2+300*n^4*r^2-1485*n^6*m^5+990*n^6*m^6-700*m^6*p^3*alpha+96
0*n^3*m^7*p+1680*n^3*m*p^3+360*n^6*m*p-555*n^5*m^4*p-675*n^2*p^4-27750
*n^3*p^3*m^2-60*n^4*m^7*alpha+300*n^8*m^2+2100*n^2*p^2*m*r-360*n^5*m*r
+3980*n^2*m^5*p^2*alpha+3780*n^5*m^2*r-1140*n^4*p^3+30*n^4*m^8+435*n^5
*m^5*alpha-180*n^7*m^2+8440*n^3*m^4*p^2-5240*n^2*m^7*r+13440*n^2*m^8*r
-180*n^5*p^2+51420*n^2*m^5*p*r+8560*n^2*m^6*r*alpha-20240*n^3*m^6*r+24
0*n^2*m^10*p-13480*n^2*p^3*m^3+30240*n^2*p^3*m^4-60*n^4*m^9+3280*m^9*n
*alpha*r-2080*n^2*m^6*p^2+8245*n^4*m^4*r+810*m^6*p^4-240*n^5*m^6+32590
*n^2*p^2*m^2*r-3600*n^2*p^2*r*alpha-6200*m^5*p^2*r-160*m^13*r-18750*m*
r^3*alpha-780*n^4*m^3*r+31600*n^2*m^3*r^2-14160*n^2*m^4*p*r+5840*n^3*m
^5*r+3000*n^3*m*r^2+540*n^5*m^7+1140*n^7*m^3+4020*n^2*p^2*m^7+40*m^12*
p^2-660*n^5*p*m^3*alpha+40*m^10*p^2-4640*n^4*p^2*m^2*alpha-1300*m^5*p^
4-1800*m^6*r^2+16930*n^3*p^2*m^4*alpha+440*m^9*p^3-3000*m^8*r^2+400*m^
7*p^3-840*m^8*p^3+5600*m^7*r^2-80*m^11*p^2-120*n^2*m^9*p-100*m^6*p^2*a
lpha^2*n-160*m^11*r+40*m^8*p^2*alpha^2+320*m^12*r+6060*n^3*p^2*r-100*m
^7*p^3*n+550*m^4*p^4-240*m^10*p^2*n+80*m^10*p^2*alpha+1340*m^5*n*p^3*a
lpha+7210*m^4*n*p*alpha^2*r+4750*n*m^2*r^2*alpha^2-31110*n^2*m^4*p*alp
ha*r-1600*m^6*p*alpha^2*r-3760*m^8*p*alpha*r+1360*m^7*n*alpha^2*r+740*
m^7*p^3*alpha+21880*m^6*p*n*r*alpha+16400*m^8*p*n*r+260*m^7*n*alpha*p^
2+6325*m*p^2*n*alpha^2*r+25100*m^4*n*r^2*alpha+18000*m^6*r^2*n-5200*m^
6*r^2*alpha-7340*m^7*r*p^2+3600*m^5*r^2*alpha-9400*m^5*r*alpha*p^2-185
20*m^5*n*p*r*alpha+25800*m^5*n*r*p^2-1750*m^4*r^2*alpha^2+1920*m^11*n*
r-320*m^11*r*alpha-160*m^9*alpha^2*r-2160*m^10*r*p+3600*m^7*p*alpha*r-
27600*m^7*p*n*r-2960*m^8*n*alpha*r+11360*m^6*n*p*r-660*m^6*p^3*n-2000*
m^8*p*r+13100*m^6*p^2*r+4160*m^9*p*r+3375*p^5*alpha+5625*p^4*r-4275*p^
5*n+320*m^10*r*alpha-3520*m^10*n*r+9400*m^4*p^2*r*alpha-160*m^8*n*p^2-
80*m^9*p^2*alpha+400*m^9*p^2*n-32760*m^4*p^2*n*r+1600*m^9*r*n-340*m^8*
p^2*n*alpha-1740*m^4*n*p^3*alpha-30000*m^5*r^2*n-19320*n*m^3*p^4-1560*
n^4*m^2*alpha*r-11700*n*m^3*r^2*alpha+11160*n*m^3*r*alpha*p^2-2300*m^3
*p^5):}{%
}
\end{mapleinput}

\end{maplegroup}
\begin{maplegroup}
\begin{mapleinput}
\mapleinline{active}{1d}{xi := evalf((-xi2+sqrt(xi2^2-4*xi3*xi1))/(2*xi3)):}{%
}
\end{mapleinput}

\end{maplegroup}
\begin{maplegroup}
\begin{mapleinput}
\mapleinline{active}{1d}{eta :=
evalf(-(4*xi*m^4-15*xi*p-21*m^2*p*alpha+10*xi*n^2+21*m^3*n*alpha-4*xi*
m^3-56*m^3*n^2-58*m*n*q-20*xi*q-23*m*p^2+4*m^6+23*p*q-4*m^7-10*xi*n*al
pha+13*xi*m*n-17*xi*m^2*n+4*xi*m^2*alpha+17*xi*m*p+25*n*p*alpha+29*m*n
^3+31*m^3*q+21*m*q*alpha-23*m*n^2*alpha+15*p^2+28*m^5*n-29*p*n^2-25*r*
alpha+81*p*m^2*n-6*n^3+35*n*r-28*m^4*p-4*m^5*alpha-52*m*n*p-35*m^2*r-2
4*n*m^4+22*n*q-26*m^2*q+30*m*r+35*m^2*n^2+26*m^3*p)/(25*m*q-20*q+20*m^
3*n+19*n*p-22*m^2*p-19*m*n^2-16*m^2*n+13*n*m*alpha+6*n^2-15*p*alpha-25
*r-4*m^3*alpha+20*m*p-4*m^5+4*m^4)):}{%
}
\end{mapleinput}

\end{maplegroup}
\begin{maplegroup}
\begin{mapleinput}
\end{mapleinput}

\end{maplegroup}
\begin{maplegroup}
\begin{mapleinput}
\mapleinline{active}{1d}{#m equal to zero, n not}{%
}
\end{mapleinput}

\end{maplegroup}
\begin{maplegroup}
\begin{mapleinput}
\end{mapleinput}

\end{maplegroup}
\begin{maplegroup}
\begin{mapleinput}
\mapleinline{active}{1d}{#w :=
evalf(sqrt(400*q^2+600*p*q-100*q*n^2+225*p^2+80*p*n^2-500*n*r-200*q*n+
60*n^3+60*n*p^2)):}{%
}
\end{mapleinput}

\end{maplegroup}
\begin{maplegroup}
\begin{mapleinput}
\mapleinline{active}{1d}{#Omega :=
evalf(45562500*p^6*q*r^2-80000000*n^2*r^2*p*q^4+36300000*n^4*r^2*p^2*q
^2+25500000*n^5*r^4*p+11520*n^9*p^3*r+291600*p^9*r*n+712800*p^5*q^3*n^
3-75937500*p^5*r^3*n-62500000*n^3*r^4*p*q-15673500*p^7*q*r*n-13410000*
p^2*n^6*r^3-140625000*p^2*r^5*n+62500000*n^2*r^5*q+25000000*q^2*n^3*r^
4-299700*p^6*n^3*q^2+25000000*n^2*r^4*q^2+25312500*n^3*r^4*p^2-6144000
*q^7*n*p+619520*q^5*n^4*p^2-51840*n^10*q*p*r+777600*p^8*n^2*r-1166400*
n^9*q*r^2+34560*n^9*p^2*r^2-32400000*p^5*q^3*r-3240000*p^6*n^3*r*q+409
6000*n^3*q^6*p-84375000*n*r^4*q*p^2+26312500*p^2*n^4*r^4+1795500*p^6*n
^3*r^2-58880*n^7*q^4*p+1350000*n^7*r^4-1600000*q^4*n^4*r^2-7680000*q^6
*n*p*r+63281250*p^4*r^4+6681600*n^3*q^5*p^2-64125000*p^4*n^3*r^3+40000
000*n^3*q^3*r^3-59375000*p*n^3*r^5-2252800*n^4*q^6+39062500*n^2*r^6+77
760*n^11*r^2-432000*n^8*p*r^3-546750*p^8*q^2+101250000*n^2*r^3*p^3*q+1
03680*n^10*r^2*p-11250000*n^5*r^4*q+1473600*p^3*n^7*r^2+207360*n^6*p^5
*r+15360*n^9*q^3*p+8869500*n^5*r^2*p^4+77760*p^7*n^4*r+5184000*q^5*p^4
+1728000*p^4*n^4*q^3-3840*p^3*n^8*q^2-36160*q^4*n^6*p^2-6696000*p^4*n^
2*q^4-650240*n^5*q^5*p-8775000*p^5*n^2*r^3-8524800*q^5*n^2*p^3-19440*p
^6*n^4*q^2+8437500*p^4*r^4*n+88080*p^4*q^3*n^5+1468800*n^8*r^2*p^2-409
6000*n^4*q^5*r-57600000*p^3*q^5*r-1280*p^4*n^7*q^2+252480*n^7*q^3*p^2-
2192000*n^5*q^4*p^2-86400000*p^4*q^4*r+(1980000*p^2*n^3*r^2*q^2-460800
*q^6*p^2-25920*n^2*p^7*r-1800000*n^4*r^2*q^2*p+168960*q^5*n^2*p^2-1944
00*p^7*r*q+10125000*r^3*p^3*q*n+1395000*p^3*r^2*n^3*q-345600*q^5*p^3-5
760*n^7*q^3*p-442800*n^5*r^2*p^3-6375000*p^2*n^2*r^4+36450*p^7*q^2-384
0*n^5*r*p^5+1350000*n^4*r^3*p^2-60480*n^3*p^6*r-4218750*p^3*r^4+216000
0*p^4*q^3*r-4050000*p^4*q^2*r^2+6480*p^6*q^2*n^2-1113750*n^2*p^5*r^2-2
45760*n^3*q^5*p+15120*p^5*q^2*n^3+307200*q^6*n*p+2880000*p^3*q^4*r+162
000*p*n^6*r^3+65280*n^5*q^4*p+1440*p^3*n^6*q^2+960*p^4*q^2*n^5-38880*n
^8*r^2*p-2970000*p^4*n^2*r^2*q-78240*p^3*n^4*q^3-3037500*p^5*r^2*q+345
600*n^2*q^4*p^3+1920*q^4*n^4*p^2-259200*q^4*n*p^4+24000*q^3*n^5*p*r-11
25000*n^3*r^4*p+162000*p^6*r^2*n-5625000*r^4*p^2*q-3840*n^6*p^2*q^3-27
360*n^3*p^4*q^3-25920*p^2*n^7*r^2+48600*p^6*q^3+2137500*p^3*n^3*r^3-14
5800*p^8*r+4687500*p*r^5*n+5062500*p^4*r^3*n-5760*n^6*p^4*r-226800*p^5
*q^3*n-174000*p^4*n^4*r^2-86400*p^2*n^5*r^2*q+25920*n^7*r*p^2*q+453600
*n^6*r^2*p*q+343200*p^4*n^4*q*r+2400000*q^3*r^2*n^2*p-10800*p^3*q^2*n^
4*r-1920000*q^4*r*p^2*n-288000*n^5*r*p^2*q^2+980100*p^6*q*r*n+384000*q
^5*n*p*r+3000000*p*n^2*r^3*q^2+3750000*n*r^4*q*p-408000*q^3*n^2*p^3*r+
17280*n^6*r*p^3*q-4500000*p^2*r^3*n^2*q-1350000*r^3*p*q*n^4+1248000*q^
3*r*p^2*n^3-1512000*p^4*q^2*n^2*r-6000000*q^3*n*p^2*r^2-192000*n^3*p*q
^4*r+122400*p^5*n^3*r*q+1093500*p^5*q^2*n*r)*w-4860000*p^7*n*r^2+15360
*n^8*r*p^4+9216000*q^7*p^2+577500*p^4*n^4*r^2*q+768000*q^3*n^5*p^2*r+1
2152000*p^3*n^5*r^2*q+32400*p^2*n^7*r^2*q+141562500*p^2*n^2*r^4*q-4040
0000*q^3*n^3*p^2*r^2+120000000*q^4*n*p^2*r^2+54337500*p^4*n^2*r^2*q^2+
173600*p^3*n^6*r*q^2-1059200*p^4*n^6*r*q-42750000*p^3*n^3*r^3*q+200700
0*p^5*n^3*r*q^2-23040*p^3*n^8*r*q-391200*p^5*n^5*r*q-82800000*r^2*p^3*
q^2*n^3-26190000*p^5*q^3*n*r-2418000*p^4*n^4*r*q^2-202500000*r^3*p^3*q
^2*n-11200000*q^4*n^3*p^2*r+26880000*q^5*n*p^2*r+126000000*q^3*n*p^3*r
^2+26000000*p*n^4*r^3*q^2-12800000*r*p*q^5*n^2+748800*n^8*q^2*p*r-2025
0000*p^4*n*q^2*r^2-540000*p*n^6*r^3*q+32000*q^4*n^5*p*r+26820000*q^3*n
^2*p^4*r+45000000*p^2*n^2*r^3*q^2-93750000*r^5*p*q*n-18750000*n^4*r^5+
11520*n^10*q^3+3645000*p^6*q^3*n-4416000*n^6*p*r*q^3-168960*n^8*q^4+12
160000*n^4*p*r*q^4+8424000*q^4*n*p^5+3888000*p^7*r*q^2-1458000*p^7*q^3
-648000*n^9*r^3-9072000*n^6*r^2*p*q^2+273600*n^7*p^2*q^2*r+121500000*p
^5*q^2*r^2+1198800*n^3*p^7*r+2816000*q^5*n^3*p*r+13824000*q^6*p^3+8100
}{%
}
\end{mapleinput}

\end{maplegroup}
\begin{maplegroup}
\begin{mapleinput}
\mapleinline{active}{1d}{
0000*p^4*q^3*r^2-54720*p^5*n^5*q^2+16706250*p^6*r^2*n^2-72900*p^8*q^2*
n+2250000*n^6*r^4+5832000*p^8*r*q-30375000*n^4*r^3*p^3+5120000*q^6*n^2
*r-2880*n^9*p^2*q^2+191250000*p^3*n^2*r^4+168750000*r^4*p^3*q+4572000*
p^5*n^4*r^2-16800000*n^5*q^3*r^2+2048000*q^7*n^2-4531200*q^6*n^2*p^2+9
26720*n^6*q^5+5120*n^8*p^2*q^3+200000*q^3*n^6*r^2-1113600*p^3*n^7*q*r+
9632000*p^3*n^5*q^2*r+29214000*n^2*r*p^5*q^2+50400000*n^4*r^2*p*q^3-22
275000*p^4*n^3*q*r^2-12546000*n^6*r^2*p^2*q+218880*n^5*r*p^6+491600*p^
4*n^6*r^2+237440*n^6*p^3*q^3-446400*q^4*n^3*p^4+614400*q^4*n^4*p^3-270
0000*n^7*p*r^3+18000000*n^5*p*r^3*q-75000000*q^2*r^4*p*n-6912000*q^6*n
*p^2+5120*p^5*n^7*r+28500000*p^2*n^4*r^3*q-80000000*q^3*n^2*p*r^3-3037
50000*r^3*p^4*q*n-35397000*p^6*q^2*r*n+15360000*q^4*n^2*p^3*r+2835000*
r^2*p^6*q*n-3392000*q^3*n^4*p^3*r+4128000*p^2*n^5*r^2*q^2-33120000*p^3
*n^3*q^3*r-1668600*p^7*n^2*r*q+91260000*p^5*n^2*r^2*q-64000*n^7*p*q^3*
r-30000000*n^3*p*r^3*q^2-15000000*n^4*r^4*q+16000000*n^3*q^4*r^2+43200
000*p^3*n*q^4*r-7511400*p^5*n^4*q*r-36000000*n^2*r^2*p^2*q^3-69120*p^2
*n^9*q*r+112500000*r^4*p^2*q^2-194400*p^7*n^2*q^2-96000*n^8*q^3*r+7560
000*n^7*q*r^3+6624000*n^7*q^2*r^2+3200000*q^5*n^2*r^2-972000*p^6*q^4-5
1840*p^4*n^6*q^2+380700*p^6*q^3*n^2-30000000*n^5*q^2*r^3+5875200*q^5*n
*p^4+2187000*p^9*r-4300000*p^3*n^5*r^3+1088000*n^6*q^4*r):}{%
}
\end{mapleinput}

\end{maplegroup}
\begin{maplegroup}
\begin{mapleinput}
\mapleinline{active}{1d}{#alpha := evalf(1/10*(-20*q-15*p+10*n^2+w)/n):}{%
}
\end{mapleinput}

\end{maplegroup}
\begin{maplegroup}
\begin{mapleinput}
\mapleinline{active}{1d}{#eta :=
evalf((-44*q*n^2+2*w*n*(1/20*(400*q^2*r-260*n^2*q*r-375*p*r^2+36*n^4*r
-27*p^3*q+195*n*p^2*r+48*n*q^2*p-4*n^3*q*p)*w/((12*n^4*q-4*n^3*p^2-88*
n^2*q^2-40*n^2*p*r+125*n*r^2+117*q*n*p^2+160*q^3-27*p^4-300*r*p*q)*n)+
1/20*(-1640*n^3*p*q^2-960*p^2*n^3*r+540*p^3*q^2+240*p*n^5*q-2925*p^3*r
*n-11550*p^2*n*r*q+3900*n^2*r*p*q+405*p^4*q-540*p^5*n-6250*n*r^3-8000*
q^3*r-80*p^3*n^4+5625*p^2*r^2+2720*q^3*n*p-6000*q^2*r*p+7500*r^2*p*q-7
20*p^2*n*q^2+sqrt(Omega)+4400*q^2*n^2*r-600*n^4*r*q+2340*p^3*n^2*q+675
0*p*n^2*r^2-540*n^4*r*p+60*p^2*n^3*q)/((12*n^4*q-4*n^3*p^2-88*n^2*q^2-
40*n^2*p*r+125*n*r^2+117*q*n*p^2+160*q^3-27*p^4-300*r*p*q)*n))+45*n*p^
2-5*p*n*w+12*n^4+8*p*n^3+54*p*q*n-100*r*q-75*p*r+5*r*w-20*n^2*r)/(-3*p
*w-40*q*n-50*n*r+8*p*n^2+12*n^3+60*p*q+45*p^2)):}{%
}
\end{mapleinput}

\end{maplegroup}
\begin{maplegroup}
\begin{mapleinput}
\mapleinline{active}{1d}{#xi :=
evalf(1/20*(400*q^2*r-260*n^2*q*r-375*p*r^2+36*n^4*r-27*p^3*q+195*n*p^
2*r+48*n*q^2*p-4*n^3*q*p)*w/((12*n^4*q-4*n^3*p^2-88*n^2*q^2-40*n^2*p*r
+125*n*r^2+117*q*n*p^2+160*q^3-27*p^4-300*r*p*q)*n)+1/20*(-1640*n^3*p*
q^2-960*p^2*n^3*r+540*p^3*q^2+240*p*n^5*q-2925*p^3*r*n-11550*p^2*n*r*q
+3900*n^2*r*p*q+405*p^4*q-540*p^5*n-6250*n*r^3-8000*q^3*r-80*p^3*n^4+5
625*p^2*r^2+2720*q^3*n*p-6000*q^2*r*p+7500*r^2*p*q-720*p^2*n*q^2+sqrt(
Omega)+4400*q^2*n^2*r-600*n^4*r*q+2340*p^3*n^2*q+6750*p*n^2*r^2-540*n^
4*r*p+60*p^2*n^3*q)/((12*n^4*q-4*n^3*p^2-88*n^2*q^2-40*n^2*p*r+125*n*r
^2+117*q*n*p^2+160*q^3-27*p^4-300*r*p*q)*n)):}{%
}
\end{mapleinput}

\end{maplegroup}
\begin{maplegroup}
\begin{mapleinput}
\mapleinline{active}{1d}{#}{%
}
\end{mapleinput}

\end{maplegroup}
\begin{maplegroup}
\begin{mapleinput}
\mapleinline{active}{1d}{#Both m and n equal to zero}{%
}
\end{mapleinput}

\end{maplegroup}
\begin{maplegroup}
\begin{mapleinput}
\mapleinline{active}{1d}{#}{%
}
\end{mapleinput}

\end{maplegroup}
\begin{maplegroup}
\begin{mapleinput}
\mapleinline{active}{1d}{#alpha := -1/5*(10*q-3*p^2+25*r)/(4*q+3*p):}{%
}
\end{mapleinput}

\end{maplegroup}
\begin{maplegroup}
\begin{mapleinput}
\mapleinline{active}{1d}{#xi :=
evalf(1/10*(7360*q^4*p+5520*p^2*q^3-4000*q^3*r+270*q^2*p^3-10000*q^2*r
^2-14100*q^2*p^2*r-12500*q*r^3+3750*q*p*r^2-10800*q*r*p^3-1161*p^5*q-8
10*p^6-1125*r^2*p^3+sqrt(1843200*q^7*p^3-172800*q^6*p^4+103680*q^5*p^6
-194400*q^4*p^7-648000*p^5*q^5-72900*p^8*q^3+28125000*q*r^5*p^3+140625
00*q^2*p^2*r^4+8437500*q*p^4*r^4-891000*q^3*p^7*r+777600*q^2*p^8*r+256
50000*q^2*p^5*r^3+93750000*q^2*r^5*p+156250000*q^2*r^6-10935*p^10*q^2+
112500000*q^3*p^2*r^4+75000000*q^3*r^4*p+22500000*q^2*r^4*p^3+9315000*
q^2*p^6*r^2+250000000*q^3*r^5-121500*q^4*p^6+67500000*q^3*r^3*p^3+1012
5000*q^3*p^4*r^2+100000000*q^4*r^4-7200000*q^5*r*p^3+2592000*q^3*r*p^6
-1674000*p^5*q^4*r+24840000*p^5*q^3*r^2+486000*p^7*r*q^2+43740*p^11*r+
6144000*q^8*p+1152000*q^7*p^2+8192000*q^9+1265625*r^4*p^6+12800000*q^7
*r^2+20480000*q^8*r+1883250*p^8*q*r^2+291600*q*r*p^9+180000000*p^2*q^4
*r^3-38400000*p^2*q^6*r-15750000*q^4*p^4*r^2-13680000*q^5*p^4*r+240000
00*q^5*p^2*r^2-80000000*q^5*p*r^3+2304000*q^6*p^3*r-128000000*q^6*p*r^
2-43520000*q^7*p*r+54000000*p^3*q^4*r^2))/((160*q^3-300*q*p*r-27*p^4)*
(4*q+3*p))):}{%
}
\end{mapleinput}

\end{maplegroup}
\begin{maplegroup}
\begin{mapleinput}
\mapleinline{active}{1d}{#eta
:=(50*q*r-15*p^2*r+125*r^2+129*p^2*q+45*p^3+92*p*q^2-80*q^2*xi-120*q*x
i*p-45*p^2*xi)/(100*q*r+80*q^2+30*p*q+9*p^3):}{%
}
\end{mapleinput}

\end{maplegroup}
\begin{maplegroup}
\begin{mapleinput}
\mapleinline{active}{1d}{#}{%
}
\end{mapleinput}

\end{maplegroup}
\begin{maplegroup}
\begin{mapleinput}
\mapleinline{active}{1d}{#Calculate d3, d2, d1, d0, a, b, c, A and B in all cases}{%
}
\end{mapleinput}

\end{maplegroup}
\begin{maplegroup}
\begin{mapleinput}
\mapleinline{active}{1d}{#}{%
}
\end{mapleinput}

\end{maplegroup}
\begin{maplegroup}
\begin{mapleinput}
\mapleinline{active}{1d}{d3 :=
evalf(9/5*m^5*q-9/5*m^4*r+24/5*m^3*r-p*alpha^3-2/25*n^3-56/25*n*p^2+36
/25*n*m^7+84/25*n*m^5+72/5*q*m^2*n-63/25*p^2*m^3-10*q*m*p-9/5*m*n^4+2/
25*p^3+6/5*n^2*alpha^2-36/25*p*m^6-96/25*p*m^4-4/25*m^9-12/25*m^7*alph
a+27/5*q*n^2*m+24/25*m^6*alpha+12/25*m^4*alpha^2+p^2+54/25*m*p^2*n-38/
5*q*n*m-5*r*alpha^2-4*q*alpha^2+2*r*m+162/25*m^3*p*alpha-87/25*m^4*p*a
lpha+11/5*q*p*alpha+14/5*m*alpha^2*p-23/5*n^2*p*alpha-3*m^2*p*alpha-4/
25*m^3*alpha^3-12/25*m^5*alpha^2+12/5*p^2*alpha-3*r*n^2-5*r*alpha-6/5*
m^2*q-3*m*q^2-314/25*n*m*alpha*p+261/25*n*m^2*alpha*p-8*n*m*alpha*q-12
/5*n^3*alpha+3*r*n-12/25*m^5*alpha-5*r*m^2+2*p*r+2/5*n*q+10*r*m*alpha+
66/25*n*m^3*alpha^2+3/5*n*m*alpha^3+19/5*n*p*alpha^2-3*m*n^2*alpha^2+4
4/5*q*n*alpha+17/5*q*m*alpha^2-42/5*q*m^2*alpha+33/25*n^2*m^2+21/5*m^3
*q*alpha-58/25*m*p^2*alpha+13/5*n*p*alpha-28/5*n^2*q+6/5*n^4+4*q^2-24/
5*q*m^4+237/25*n^2*m^4-186/25*n^3*m^2+252/25*n^2*m*p+279/25*m^2*p*n-38
/5*m*n*r-23/5*m*p^2-69/25*p*n^2-12/5*m*p*n+27/5*r*m^2*n+27/5*m^2*q*p+2
1/5*m^3*q-12/5*m*p*r-12/5*q*p*n+23/5*q*p+69/25*m*n^3+3*q*r-162/25*m^3*
n^2-108/25*m^5*n^2+123/25*m^3*n^3+102/25*m^5*p-78/5*m^3*p*n+153/25*m^2
*p^2-96/25*m^6*n-12/25*m^7+4/25*m^6+21/5*m*alpha*n^3+234/25*m^2*alpha*
n^2+63/25*n*m^3*alpha-36/5*n*m^3*q-63/25*m*alpha*n^2-183/25*m^3*alpha*
n^2-189/25*m^2*n^2*p+7*r*n*alpha-29/5*r*m^2*alpha+21/5*q*m*alpha+171/2
5*n*m^4*p-6*n*m^4*alpha+87/25*n*m^5*alpha-54/25*n*m^2*alpha^2-76/25*m^
2*alpha^2*p+12/25*m^8-24/25*m^4*n+6/5*m^3*p+9/5*p*n^3):}{%
}
\end{mapleinput}

\end{maplegroup}
\begin{maplegroup}
\begin{mapleinput}
\mapleinline{active}{1d}{d2 :=
evalf(-18/5*q*m^2*eta+246/25*m^5*q+12/25*eta*m^8-138/25*m^6*q-54/5*m^4
*r+28/5*m^3*r+694/25*m*n^2*p*alpha-308/25*m*q*p*alpha+162/25*m^3*eta*p
*alpha+234/25*eta*m^2*alpha*n^2-10*eta*q*m*p+72/5*eta*q*m^2*n+24/5*m^7
*p-24/25*m^7*eta-28/5*n*p^2-24/5*n*m^8+276/25*m^4*p^2+216/25*n*m^7+6/5
*q*n*eta+292/25*q*m^2*n-468/25*p^2*m^3-48/5*q*m*p-198/25*m*n^4-12/5*p^
3-6/5*n^3*alpha^2-228/25*p*m^6+12/25*m^10+6*m^5*r-6/5*n^5+3*p^2*alpha^
2-24/25*m^9-96/25*eta*m^6*n+153/25*eta*m^2*p^2-78/5*eta*m^3*p*n+102/25
*eta*m^5*p+252/25*eta*n^2*m*p-186/25*eta*n^3*m^2-24/5*eta*q*m^4+237/25
*eta*n^2*m^4-28/5*eta*n^2*q-24/25*m^7*alpha+686/25*q*n^2*m+5*r^2+656/2
5*m*p^2*n+3*eta*p^2+12/5*p^2*xi+88/25*m^3*alpha^2*p-24/25*m^5*eta*alph
a-186/25*m^4*p*alpha+46/5*q*p*alpha-48/5*n^2*q*alpha-12/5*n^3*eta*alph
a-38/5*n^2*p*alpha-10*r*alpha*eta+24/25*m^4*alpha*xi+24/25*m^8*alpha+1
2/25*m^6*alpha^2-12/25*m^7*xi+24/25*m^6*xi-42/5*r*n^2-5*r*xi-10*m*q^2+
14/5*p*m*alpha^2*eta-152/25*m^2*alpha*p*xi+34/5*q*m*alpha*xi+44/5*q*n*
eta*alpha-42/5*q*m^2*eta*alpha-38/5*r*m*n*eta+10*r*m*alpha*eta+22/5*n*
q*alpha^2+42/5*q*m*alpha*eta+572/25*n*m^2*alpha*p-76/5*q*m*n*eta-314/2
5*n*m*xi*p+126/25*n*m^3*alpha*eta-108/25*n*m^2*alpha*xi+261/25*n*m^2*x
i*p-452/25*n*m*alpha*q-126/25*m*alpha*n^2*eta-116/5*r*m*alpha*n+558/25
*n*m^2*eta*p+12/5*n^4*alpha-76/5*r*m*q-12/5*n^3*xi+8*p*r+8/5*q^2*alpha
+51/5*m^2*q^2-26/5*q^2*n+28/5*m*alpha*p*xi+26/5*n*eta*p*alpha-6*m^2*et
a*p*alpha-816/25*n*m^3*alpha*p+638/25*n*m^2*alpha*q-8*n*m*xi*q-54/25*n
*eta*m^2*alpha^2-6*n*eta*m^4*alpha+9/5*n*m*alpha^2*xi+132/25*n*m^3*alp
ha*xi-42/5*n*m*p*alpha^2+38/5*n*p*alpha*xi-6*m*n^2*alpha*xi-314/25*n*e
ta*m*alpha*p+6*r*n*eta-10*r*m^2*eta+10*r*m*xi-78/25*n*m^4*alpha^2+12/5
*n^2*alpha*xi+6/5*n^2*eta*alpha^2+129/25*m^2*n^2*alpha^2-10*p^2*n*alph
a+12/5*p^2*alpha*eta-528/25*m^3*q*p+268/25*m^2*p^2*alpha-42/5*q*m^2*xi
+44/5*q*n*xi+24/5*r*m^3*eta+6/5*eta*n^4-3*m^2*xi*p+4*eta*q^2-63/25*m*x
i*n^2-87/25*m^4*p*xi+162/25*m^3*p*xi+11/5*q*p*xi+216/25*m^3*q*alpha-8*
m*p^2*alpha-198/25*n*m^6*alpha+516/25*m^4*alpha*n^2-444/25*m^2*alpha*n
^3+21/5*m^3*q*xi-58/25*m*p^2*xi+13/5*n*p*xi-23/5*n^2*p*xi+24/25*m^6*et
a*alpha+198/25*m^5*p*alpha+12/25*m^4*alpha^2*eta-3*p*alpha^2*xi-12/25*
m^5*xi+63/25*m^3*xi*n-6/25*n^3*eta+171/25*p^2*n^2+99/25*n^2*m^2*eta-32
/25*n^2*q+6/25*n^4+4/5*q^2-108/25*q*m^4+234/25*n^2*m^4-168/25*n^3*m^2+
312/25*n^2*m*p-56/5*m*n*r+52/25*m*p^3+84/5*m^2*p*r+96/5*m*r*n^2+152/5*
r*m^2*n-58/5*p*n*r+748/25*m^2*q*p+18/5*m^3*eta*p-72/25*m^4*eta*n+654/2
5*q*n*m*p-104/5*m*p*r-476/25*q*p*n+12/25*m^6*eta+10*q*r+6*q*n^3-36/5*n
*eta*m*p-642/25*m^5*n^2+696/25*m^3*n^3+108/25*m^5*p-432/25*m^3*p*n+38/
5*m^2*p^2-96/25*m^6*n-52/25*q*p^2+321/25*n^4*m^2-24*m^3*r*n-88/25*m^2*
alpha^2*q-216/25*m^4*q*alpha+186/25*m*alpha*n^3+228/5*m^3*p*n^2-138/25
*n^2*eta*p-792/25*q*m^2*n^2-183/25*m^3*xi*n^2+21/5*m*xi*n^3-324/25*m^3
*n^2*eta-492/25*n^3*m*p+87/25*n*m^5*xi+678/25*n*q*m^4-192/5*n*m^3*q+23
4/25*m^2*xi*n^2-72/5*m^3*alpha*n^2-1296/25*m^2*n^2*p+138/25*m*n^3*eta+
56/5*r*m^3*alpha+6*r*m*eta+10*r*n*alpha+10*r*p*alpha+4*r*m*alpha^2-48/
5*r*m^2*alpha+7*r*n*xi-29/5*r*m^2*xi+1122/25*n*m^4*p+168/25*n*m^5*eta-
6*n*m^4*xi+174/25*n*m^5*alpha+42/5*q*m^3*eta+46/5*q*p*eta-576/25*n*m^2
*p^2-702/25*n*m^5*p-192/25*p*m^4*eta-24/25*m^5*alpha*xi-46/5*m*p^2*eta
-12/25*m^3*alpha^2*xi+2*p*r*eta-8*q*alpha*xi-4*q*eta*alpha^2-10*r*alph
a*xi-56/25*p^2*n*eta+84/5*m^6*n^2+12/25*m^8-606/25*n^3*m^4+21/5*m*xi*q
+198/25*p*n^3):}{%
}
\end{mapleinput}

\end{maplegroup}
\begin{maplegroup}
\begin{mapleinput}
\mapleinline{active}{1d}{d1 :=
evalf(-18/5*q*m^2*eta^2+6/5*q*n*eta^2-12/25*m^7*eta^2+31/5*p^3*n-36/5*
n*eta^2*m*p+24/25*eta*m^8-126/25*m^6*q-31/5*p^2*r+99/25*n^2*m^2*eta^2+
18/5*m^3*eta^2*p-72/25*m^4*eta^2*n-3*p^3*alpha+87/5*n*p^2*m*alpha+367/
25*m*n^2*q*alpha+152/5*r*n*m^2*eta+579/25*m^4*p*n*alpha-176/25*m^2*alp
ha*q*xi-72/5*eta*m^3*alpha*n^2-1296/25*eta*m^2*n^2*p-63/25*m*alpha*n^2
*eta^2+202/5*r*p*m*n+279/25*n*m^2*eta^2*p-126/25*eta*m*xi*n^2-6*eta*m^
2*xi*p-282/25*q^2*m^3+141/25*m^7*q-24/25*m^9*eta-96/5*eta*q*m*p+584/25
*eta*q*m^2*n+126/25*m^7*p-132/25*p*m^8-24/5*n*m^8-186/25*p^3*m^2+327/2
5*m^4*p^2-15*p^2*m^5-129/25*p*n^4-12/25*m^11+12/25*m^10+12/25*m^6*eta^
2-33/5*m^6*r+6*m^5*r-11*m*r^2-6/25*n^5+42/5*eta*m*xi*q-192/25*eta*m^6*
n+76/5*eta*m^2*p^2-864/25*eta*m^3*p*n+216/25*eta*m^5*p+624/25*eta*n^2*
m*p-336/25*eta*n^3*m^2-216/25*eta*q*m^4+468/25*eta*n^2*m^4-64/25*eta*n
^2*q+129/25*m*n^5+126/25*eta*m^3*xi*n+6/5*n^2*xi^2+5*r^2-4*q*xi^2-12/5
*p^3*eta-5*r*xi^2+3*p^2*eta^2-69/5*n^2*m^5*alpha-37/25*m*q^2*alpha-123
/25*m*n^4*alpha-822/25*m^5*n*q-12/25*m^9*alpha-452/25*n*eta*m*alpha*q-
24/25*eta*m^5*xi-24/25*m^7*xi+927/25*m^5*n^3-528/25*m^7*n^2+44/5*n*q*a
lpha*xi+13/5*p*n*eta^2*alpha-3*p*m^2*eta^2*alpha+46/5*p*q*alpha*eta+24
/25*m^4*alpha*eta*xi-104/5*r*m*p*eta-112/5*r*m*n*eta-762/25*n^2*m^2*p*
alpha-186/25*m^4*eta*p*alpha+10*r*m*xi*eta+10*r*n*eta*alpha+8*r*m*alph
a*xi-192/5*n*eta*m^3*q+656/25*n*eta*m*p^2+21/5*q*m*alpha*eta^2-38/5*q*
m*n*eta^2+572/25*n*m^2*xi*p+63/25*n*m^3*alpha*eta^2+186/25*eta*m*alpha
*n^3+686/25*eta*m*n^2*q+79/5*r*m^2*n*alpha-116/5*r*m*xi*n+132/25*m^9*n
-669/25*m^3*n^4-52/5*r*m*q+12/25*m^4*xi^2+14/5*m*p*xi^2-54/25*m^2*n*xi
^2+147/25*m^2*q^2-34/25*q^2*n+28/5*p*m*alpha*xi*eta-8*q*eta*alpha*xi-8
*p^2*m*alpha*eta-10*p*n*q*alpha-38/5*p*n^2*eta*alpha+259/25*q*m^2*p*al
pha-444/25*q*m^3*n*alpha+216/25*q*m^3*eta*alpha+748/25*p*q*m^2*eta+26/
5*n*eta*p*xi+176/25*m^3*xi*p*alpha+174/25*n*eta*m^5*alpha-452/25*n*m*x
i*q+1122/25*n*eta*m^4*p-1404/25*p*q*m^2*n-42/5*p*r*m*alpha+9/5*n*m*alp
ha*xi^2-156/25*n*m^4*alpha*xi+12/5*n^2*eta*alpha*xi+258/25*m^2*n^2*alp
ha*xi-48/5*r*m^2*eta*alpha-476/25*n*eta*q*p-108/25*n*eta*m^2*alpha*xi-
84/5*n*m*p*alpha*xi-42/5*n^2*r*eta+6*p^2*alpha*xi+10*r*q*eta-5*r*m^2*e
ta^2+3*r*n*eta^2-132/5*p*r*m^3+19/5*n*p*xi^2-12/5*n^3*alpha*xi-10*p^2*
n*xi+12/5*p^2*xi*eta-3*p*alpha*xi^2-198/5*m^2*r*n^2-54/5*m^4*r*eta+246
/25*m^5*eta*q-504/25*m^3*q*p+217/25*q*p^2*m+24/25*m^6*alpha*xi+111/25*
m^7*n*alpha-222/25*m^3*p^2*alpha-111/25*m^6*p*alpha-24/25*m^7*eta*alph
a+56/5*r*m^3*eta+12/25*eta*n^4+8/5*eta*q^2-186/25*m^4*p*xi+46/5*q*p*xi
+216/25*m^3*q*xi-8*m*p^2*xi-38/5*n^2*p*xi-5*r*alpha*eta^2+5*p*n^3*alph
a-10*r*xi*eta-6/25*n^3*eta^2+47/5*p^2*n^2+26/5*m*p^3+94/5*m^2*p*r+94/5
*m*r*n^2-84/5*p*n*r+748/25*q*n*m*p+34/25*q*n^3-26/5*q*p^2-3*m*xi^2*n^2
-31/25*p*q^2+267/25*n^4*m^2-76/25*m^2*xi^2*p+27/5*r*n^3+66/25*m^3*xi^2
*n-24*m^3*r*n+1314/25*n^2*q*m^3-27/5*m^4*r*alpha+33*m^4*r*n-516/25*n^3
*q*m-846/25*n^2*p^2*m+1218/25*n*p^2*m^3+252/5*m^3*p*n^2-198/25*eta*m*n
^4+69/25*m*n^3*eta^2-69/25*n^2*eta^2*p-132/5*q*m^2*n^2-72/5*m^3*xi*n^2
-162/25*m^3*n^2*eta^2+186/25*m*xi*n^3-504/25*n^3*m*p+301/25*n*q^2*m+17
4/25*n*m^5*xi+606/25*n*q*m^4+84/25*n*m^5*eta^2+696/25*eta*m^3*n^3-642/
25*eta*m^5*n^2+198/25*eta*n^3*p+381/25*q*p*n^2-54/5*q*n*r+109/5*q*m^2*
r+6*r*m*eta^2+56/5*r*m^3*xi+10*r*n*xi+10*r*p*xi+r*q*alpha-48/5*r*m^2*x
i+21/5*q*m^3*eta^2+23/5*q*p*eta^2-5*r*n^2*alpha-804/25*n*m^2*p^2-756/2
5*n*m^5*p+216/25*n*eta*m^7+183/5*p*m^6*n-468/25*p^2*m^3*eta-96/25*p*m^
4*eta^2-228/25*p*m^6*eta+657/25*p*q*m^4-1959/25*p*n^2*m^4+1362/25*p*n^
3*m^2+402/25*n^3*m^3*alpha+111/25*q*m^5*alpha-12/25*m^5*eta^2*alpha-23
/5*m*p^2*eta^2-10*m*q^2*eta-12/25*m^3*alpha*xi^2+16*p*r*eta-56/5*p^2*n
*eta+572/25*n*eta*m^2*alpha*p+12/5*n^4*xi+8/5*q^2*xi-12/25*m^5*xi^2+41
}{%
}
\end{mapleinput}

\end{maplegroup}
\begin{maplegroup}
\begin{mapleinput}
\mapleinline{active}{1d}{
7/25*m^6*n^2-582/25*n^3*m^4-12/5*n^3*eta*xi-308/25*q*m*p*xi+638/25*q*m
^2*n*xi+44/5*q*n*eta*xi-42/5*q*m^2*eta*xi+24/25*m^8*xi+17/5*m*xi^2*q-1
98/25*m^6*n*xi+268/25*m^2*p^2*xi-816/25*m^3*p*n*xi+198/25*m^5*p*xi-314
/25*n*eta*m*p*xi+24/25*m^6*eta*xi-6*m^4*eta*n*xi+162/25*m^3*eta*p*xi+6
94/25*n^2*m*p*xi-444/25*n^3*m^2*xi+516/25*n^2*m^4*xi-216/25*q*m^4*xi-4
8/5*n^2*q*xi+234/25*n^2*m^2*eta*xi):}{%
}
\end{mapleinput}

\end{maplegroup}
\begin{maplegroup}
\begin{mapleinput}
\mapleinline{active}{1d}{d0 :=
evalf(14/25*q^2*n^2-528/25*n^2*q*m*p-12/25*n^4*q+408/25*n^3*q*m^2-72/5
*m^5*r*n-24*m^3*r*n*eta+6*m^5*r*eta-27/5*m^4*r*xi+24*m^3*r*n^2+p^4-708
/25*n*p^2*m^4+87/5*n*p^2*m*xi+47/5*n^2*p^2*eta-132/5*n^2*q*m^4+12/5*m^
7*r-5*r*n^2*xi-48/5*m^3*r*q-24/5*n^3*p^2-37/25*m*q^2*xi-123/25*m*n^4*x
i+579/25*m^4*p*n*xi-186/25*m^4*eta*p*xi+259/25*q*m^2*p*xi-444/25*q*m^3
*n*xi+216/25*q*m^3*eta*xi+44/5*r*p*n^2-132/5*eta*q*m^2*n^2+186/25*eta*
m*xi*n^3-63/25*m*xi*n^2*eta^2+234/25*m^4*eta^2*n^2+33/25*m^2*eta^3*n^2
-168/25*n^3*m^2*eta^2+312/25*n^2*eta^2*m*p-32/25*q*n^2*eta^2-432/25*n*
m^3*eta^2*p-54/25*n*eta*m^2*xi^2-24/5*n*eta*m^8+6*m^2*r^2+63/25*n*m^3*
xi*eta^2+292/25*n*q*m^2*eta^2+6/25*n^4*eta^2-6/25*eta*n^5-582/25*eta*n
^3*m^4+34/25*eta*n^3*q-72/5*eta*m^3*xi*n^2+267/25*eta*n^4*m^2+252/5*et
a*m^3*p*n^2-504/25*eta*n^3*m*p+56/5*q*m*n*r+417/25*eta*m^6*n^2-2/25*n^
3*eta^3-4*n*r^2+4/25*q^3-48/5*n*p^3*m+804/25*n^2*p^2*m^2+2*r*m*eta^3+3
67/25*m*n^2*q*xi+r*q*xi+10*r*n*eta*xi-48/5*r*m^2*eta*xi+79/5*r*m^2*n*x
i+327/25*m^4*p^2*eta-168/25*m^2*p^2*q-8/5*q*p*r-756/25*n*eta*m^5*p+606
/25*n*eta*q*m^4+174/25*n*eta*m^5*xi-108/25*q*m^4*eta^2-452/25*n*eta*m*
xi*q-96/25*n*m^6*eta^2+21/5*q*m*xi*eta^2-804/25*n*eta*m^2*p^2-24/25*n*
m^4*eta^3+136/25*m^3*p^3+2/5*q*n*eta^3-6/5*q*m^2*eta^3-48/5*q*m*p*eta^
2+162/25*m^6*p^2+4/5*q^2*eta^2+748/25*n*eta*q*m*p+572/25*n*eta*m^2*xi*
p-24/25*m^7*eta*xi-12/25*m^5*eta^2*xi+28/5*q*n*p^2-84/5*p*r*n*eta-144/
5*p*r*m^2*n-42/5*p*r*m*xi+12*p*r*m^4+28/25*q^2*m*p-42/5*n*m*p*xi^2-78/
25*n*m^4*xi^2+3/5*n*m*xi^3+6/5*n^2*eta*xi^2+129/25*m^2*n^2*xi^2-6/5*n^
3*xi^2+5*r^2*eta+4*r*m*xi^2+4/25*m^6*eta^3+12/25*m^8*eta^2+147/25*m^2*
q^2*eta-34/25*q^2*n*eta+22/5*n*q*xi^2-4*q*eta*xi^2+94/5*p*r*m^2*eta-38
/5*p*n^2*eta*xi-10*p*n*q*xi+108/25*m^5*eta^2*p+126/25*m^7*eta*p-126/25
*m^6*eta*q+12/25*m^10*eta-504/25*m^3*eta*q*p-3*p*m^2*eta^2*xi+13/5*p*n
*eta^2*xi+46/5*p*q*xi*eta+14/5*p*m*xi^2*eta+102/25*m^4*q^2-48/25*m^8*q
+4/25*m^12+111/25*m^7*n*xi-p*xi^3-52/5*r*m*q*eta-56/5*r*m*n*eta^2+94/5
*r*m*n^2*eta+28/5*r*m^3*eta^2-48/5*m*n^3*r+2/25*n^6+3*p^2*xi^2-3*p^3*x
i+6/5*m^3*p*eta^3-12/5*m*p*n*eta^3+28/5*m*p^2*r+8*p*r*eta^2+816/25*m^3
*p*q*n+48/25*m^9*p+26/5*p^3*m*eta-26/5*p^2*q*eta+38/5*m^2*p^2*eta^2-38
4/25*m^7*p*n-264/25*m^5*p*q+204/5*m^5*p*n^2-1008/25*m^3*p*n^3+252/25*m
*p*n^4-28/5*p^2*n*eta^2-8*p^2*m*xi*eta-762/25*n^2*m^2*p*xi+402/25*n^3*
m^3*xi+216/25*m^8*n^2-48/25*m^10*n-448/25*m^6*n^3+417/25*m^4*n^4-168/2
5*m^2*n*q^2+324/25*m^6*n*q-132/25*m^2*n^5-69/5*n^2*m^5*xi+111/25*q*m^5
*xi+12/25*m^4*xi^2*eta+88/25*m^3*xi^2*p-88/25*m^2*xi^2*q-222/25*m^3*p^
2*xi-111/25*m^6*p*xi-4/25*m^3*xi^3+12/25*m^6*xi^2-12/25*m^9*xi-5*r*xi*
eta^2+5*p*n^3*xi+p^2*eta^3):}{%
}
\end{mapleinput}

\end{maplegroup}
\begin{maplegroup}
\begin{mapleinput}
\mapleinline{active}{1d}{d :=
evalf(1/6*(36*d1*d2*d3-108*d0*d3^2-8*d2^3+12*sqrt(3)*sqrt(4*d1^3*d3-d1
^2*d2^2-18*d1*d2*d3*d0+27*d0^2*d3^2+4*d0*d2^3)*d3)^(1/3)/d3-2/3*(3*d1*
d3-d2^2)/(d3*(36*d1*d2*d3-108*d0*d3^2-8*d2^3+12*sqrt(3)*sqrt(4*d1^3*d3
-d1^2*d2^2-18*d1*d2*d3*d0+27*d0^2*d3^2+4*d0*d2^3)*d3)^(1/3))-1/3*d2/d3
):}{%
}
\end{mapleinput}

\end{maplegroup}
\begin{maplegroup}
\begin{mapleinput}
\mapleinline{active}{1d}{b := evalf(alpha*d+xi):}{%
}
\end{mapleinput}

\end{maplegroup}
\begin{maplegroup}
\begin{mapleinput}
\mapleinline{active}{1d}{c := evalf(d+eta):}{%
}
\end{mapleinput}

\end{maplegroup}
\begin{maplegroup}
\begin{mapleinput}
\mapleinline{active}{1d}{a :=
evalf(1/5*d*m^3+1/5*b*m-3/5*d*m*n-2/5*n^2+4/5*q-1/5*c*m^2+2/5*c*n-4/5*
m*p+4/5*m^2*n+3/5*d*p-1/5*m^4):}{%
}
\end{mapleinput}

\end{maplegroup}
\begin{maplegroup}
\begin{mapleinput}
\mapleinline{active}{1d}{A :=
evalf(m^2*r*b^3+6*c^2*r*b*d*p-10*c^2*r*b*a+r*d^3*q^2-2*r*b^3*n+15*a^2*
b*r-10*q^2*a*d*m*n-q^2*c*n*d*p-4*d^2*q^3*c-9*c*r*p*b^2+18*q^2*a^2-8*a*
q^3-16*q*a^3-6*c*r*b^2*d*n-4*c^2*r*b^2*m+3*c*r*b*d^2*n^2+6*c*r*b*p^2-5
*c^3*r^2-6*c^2*r*d*p^2+3*c^3*r*b*n-5*m*q^2*b^2*d-8*m*r*a*c^2*n+4*q^2*a
*n^2-8*m*q*b^2*r-2*m^2*q*a*b^2+d^4*q^3+3*r^2*d^2*n^2-2*d^3*n*q^2*b-5*d
^3*r^2*b+4*m*r*a*c*n^2-d*r*b^3*m+r*d*m^2*n*b^2-5*r*d*p*b^2*m-r*b*d*p*n
^2-2*r*b*a*m^2*n-2*c^3*q^2*n-4*r*d*q*a*c+r*b*d^2*p*m*n+2*r*b*a*n^2+2*r
*a*b^2*m+2*r*m^2*p*b^2-r*n^2*b^2*m+r*p*b^2*n+2*r^2*m*p*c+r*b*n*p^2+3*b
^2*r*d^2*p+4*d*r*a*c*m*p-d^3*r*b*n*p-d^2*r*b^2*m*n+b^4*q-2*b^2*q*d^2*m
*p-3*m*q^3*b+d*m^2*q*b^3+2*d*r*c^3*m*p-8*a*n*r^2-2*b^2*q*c*n^2+3*b^2*q
*c*d*p+b^2*q*c^2*n-2*a*d^3*p^3-2*c^2*r*p*b*m+2*c*r*d^3*p^2+10*a*c*r^2-
12*q*a^2*n^2-7*p*b*r^2-2*c*q^3*n+13*d^2*r^2*b*m+c^2*q*b*d*m*p+4*c^2*r*
m*p^2+4*c^3*r*a*m-2*c^2*r*p*d^2*n+4*q*n*r^2-b*q*c^3*p+b*q*d^3*p^2-c^2*
q*p*b*m^2-4*m^2*r^2*c^2+5*d^3*r*b*m*q-8*c^2*q^2*a+4*d*r*b^2*c*m^2-4*d^
2*r*b*m*c*p-m^3*q*b^3+m^2*q^2*c^3+3*m^2*q^2*b^2+4*c*q^2*b^2-2*m^2*r*c^
3*p-4*d^3*r^2*c*m-3*d*r*b*c^2*m*n-9*d^2*r^2*c*n-2*d^2*r*m*c*p^2+3*m^2*
r*b*c^2*n+8*r^2*c^2*n+10*d^2*r^2*a-6*a*p^3*b+4*d*r*p*b*c*m^2-4*m^3*r*b
^2*c-3*n*q^2*b^2+q^2*c^2*n^2+6*a*p^2*b^2+3*d^4*r^2*n-2*c^4*r*p+b^2*q*n
^3-c*q^2*d^3*p+6*m^2*q*c*a^2+c^4*q^2+d*r*c^3*q+12*a*m^2*r^2+r^2*d^2*q+
q^2*c*p^2+7*q^2*b*r+12*a^2*c*n*m*p+3*a^2*c*n*d*p-15*a^2*c*m*d*q-3*a^2*
c*n*b*m+3*a^2*d*m^2*n*b-9*a^2*d^2*p*m*n-3*a^2*m^3*n*b+9*a^2*d*p*m^2*n+
12*a^2*n*p^2-12*a^2*p*r+3*a^2*d^2*n^3+9*a^2*d^2*p^2+6*a^2*d*p*c*m^2+24
*a^2*m*n*r-3*a^2*d*m*n^3-9*a^2*d^2*n*q-12*a^2*m*p*n^2-6*a^2*m^3*p*c+9*
a^2*n^2*b*m+3*a^2*d*p*n^2-15*a^2*d*p^2*m+9*a^2*d^2*m^2*q-6*a^2*c^2*m*p
-3*a^2*c*n^2*d*m+4*a^3*m^4+8*a^3*n^2+3*a^2*n^4-6*q^2*a*c*m^2+4*q^2*c*n
*a+12*a^3*d*m*n+3*m^2*p*a^2*b+21*a^2*d*m^2*r+3*a^2*c*n^2*m^2-6*a^2*b*d
*n^2-15*a^2*p*b*n+9*a^2*b*c*p+12*a^2*b*d*q-18*a^2*c*m*r-21*a^2*d*n*r+5
*a^4-2*q^2*d*p*a+8*q^2*a*m*p+4*a*c*m*d*q^2+6*a*d^2*m^2*q^2+q^2*d*n*r+3
*q^2*p*b*n-5*q^2*b*c*p+3*c*m*d*q^3-16*a*b*d*q^2+4*a*d^2*n*q^2+d^2*n*q^
3+3*a^2*c^2*n^2+4*a^3*c*m^2-4*a^3*b*m-8*c*n*a^3-12*d*p*a^3-4*a^3*d*m^3
+16*a^3*m*p-16*a^3*m^2*n+3*a^2*b^2*n-6*a^2*c*n^3-12*a^2*m^3*r-9*a^2*c*
p^2-4*q^2*p*r+q^4+4*b*d*q^3-d*p*q^3-3*d*p*a^2*b*m+4*c*q^2*b*d*n-3*d*q^
2*b*c*m^2-d*q^2*c^2*m*n+c^2*q^2*b*m+d^2*p*b*q^2+15*d*r*c*a^2-9*d^2*r*m
*a^2+2*c^2*q^3-d^3*q^3*m+q^2*c*n*b*m+d^2*p*c*m*q^2-m*p*b*d*q^2+12*a*m*
p*b*r-3*r*c*d*q^2-r*d^2*m*q^2+3*c^2*p*d*q^2+2*a*c^3*p^2-2*a*m^2*p*c*n*
b+2*a*m*p*b*d*n^2-2*a*m*p^2*c*n*d+16*a*c^2*p*r+2*a*c*p^2*n^2+2*a*m^2*p
^2*c^2-16*q*a*c*n*d*p+2*q*m*p*b*r+10*a*m*q^2*b+4*m^3*r*a*c^2+6*q*a*d*p
*n^2+2*q*a*d*p^2*m+2*q*a*n^2*b*m-16*q*a*c*n*b*m-2*q*a*d*m^2*n*b-6*q*a*
d^2*p*m*n+8*q*a*c*n*m*p+4*q*a*p*b*n+4*q*a*b*c*p+8*q*a*c*m*r+12*q*a*d*n
*r-16*q*a*m*n*r+4*q*a*c*n^2*d*m+2*q*a*d*p*c*m^2+10*q*a*d*m^2*r-4*q*a*c
*n^3+4*q*a*c*p^2-8*q*a*n*p^2+16*q*a*p*r-2*q*a*d^2*p^2+6*q*c*n*a^2+15*q
*d*p*a^2-9*q*a^2*d*m^3-24*q*a^2*m*p+12*q*a^2*m^2*n+4*q*a*b^2*n-3*q*a^2
*b*m-10*q*m^2*p*a*b+6*q*a^2*d*m*n-4*q*c^2*p*r+8*q*a*c^2*n^2+6*c^2*q*a^
2+6*a^2*m^2*p^2-3*q*m*p*b^2*n-16*q*d^2*r*m*a-q*b*p*c*n^2+2*q*b*p*c^2*n
-q*d^2*r*n*c*m-22*q*a*b*r+6*a*d*p^3*c-2*a*d*p^3*n+2*a*d*p^2*r+11*q*d^2
*r*c*p-2*q*d*r*c^2*n+d^3*r*c*q*n+q*d*r*c*n^2+3*q*d^3*r*m*p-5*q*d^2*r*b
*n-5*q*d^2*r*b*m^2-q*c*p*b^2*m+4*q*c*p*n*r+13*q*d*r*b^2-q*d^2*m*p^2*b+
3*q*p^2*b^2-6*q*c*r^2-3*q*b^3*p+2*q*d*m^2*p*b^2+10*q*r*b*d*m*n-2*a*b^3
*p-q*d*m*r^2+3*q*d*p^2*r+20*q*d*p*a*b*m+q*d*p*c*n*b*m+2*q*m*p^2*b*c-8*
}{%
}
\end{mapleinput}

\end{maplegroup}
\begin{maplegroup}
\begin{mapleinput}
\mapleinline{active}{1d}{
q*d*p*c*m*r-3*q*d^2*p*n*r+q*d*p*b^2*n+5*d*r^3-4*m*r^3+2*a*p^4+4*a*d*p*
b^2*n-6*a*d^2*p*m^2*r+2*a*d^2*p*n*r+6*a*m*p^2*b*n+2*a*d*p*c*n*b*m-8*a*
m*p^2*r+2*a*m*p^3*d^2-3*q*d*p^2*b*c+2*a*m*p^2*b*c+q*d*p^2*b*n-2*a*d*p^
2*b*n-2*a*d^2*p*b*n^2-10*a*d^2*p*b*q-6*a*d*p^2*b*c+8*a*c*q*b*d*n-4*a*c
*q*d^2*n^2+6*a*c^2*q*b*m-6*a*m*p*b^2*n+2*a*d^2*p^2*c*n-2*a*d*p^2*c^2*m
+6*a*d^3*p*n*q+4*a*b*p*c*n^2-2*a*b*p*c^2*n-2*a*b*p*n^3+4*a*m^2*p*c*r-8
*a*c*q*b^2+14*a*d^2*r*b*n-8*a*d^2*r*b*m^2-14*a*d^2*r*c*p-6*a*d*r*c^2*n
+14*a*d*r*b*c*m+12*a*d*r*c*n^2+6*a*d^3*r*m*p+10*a*d^2*r*n*c*m-4*a*c^3*
q*n-4*a*m^2*q*c^2*n+6*a*m^3*q*b*c+2*a*c^2*p*d*q+2*a*c*p*b^2*m-24*a*c*p
*n*r-6*a*d*q*b*c*m^2+4*a*d*q*c^2*m*n+2*a*m^3*p*b^2-4*a*m^2*p^2*b*d-6*a
*d*r*n^3-6*a*d^3*r*n^2-10*a*d*r*b^2+6*a*r*d^2*n^2*m-10*a*r*d*n*c*m^2-2
0*a*r*b*d*m*n+8*a*r*b*d*m^3-2*a*d*m^2*p*b^2+4*a*d^2*m*p^2*b-4*a*m*p^3*
c+8*a*r*p*n^2-6*a*d^3*q^2*m+8*a*d^2*q^2*c-22*a*d*m*r^2-4*a*c^2*p^2*n-2
*a*d^2*p*c*m*q+b*d*n*r^2-5*p*d*n*r^2-2*r*p*c^2*n^2+5*b*m*n*r^2-8*b*d*m
^2*r^2-2*b*c*m*r^2+4*c*p*d*r^2+2*r*c^2*p*d*n*m-5*r*b*q*n^2+4*r*c^3*p*n
-4*d*r*a*c^2*m^2+2*r*b*c*n*q+2*r*b*d*p*a-3*r*b*c*n^2*d*m+3*r*b*c*n^3-6
*r*b*c^2*n^2+2*r*m*q^2*c-8*r*p*c*n*b*m+5*b^2*r^2-5*d*r^2*c*b-10*r*p*b*
d*q+2*r*c*n*d*p^2+4*d^2*r^2*c*m^2-3*d^3*r^2*m*n+3*d*r^2*c^2*m+5*r*b*c*
n*d*p-6*r*b*a*c*m^2+8*r*b*c*n*a+r*d*q*c^2*m^2-2*c*r*p^3+3*r*c^2*q*b+11
*r*c*n*b^2*m+4*r*p*a*d*m*n+2*r^2*d*n*c*m+2*r^2*d^2*m*p-3*r^2*c*n^2+2*p
^2*r^2-2*d^3*r^2*p+2*d^2*q^2*m*b*n+3*m*q*b^3*n+2*m^2*q*b*c*r-2*d*q^2*b
*n^2+5*c*r*b^3-m*q*b^2*d*n^2+m^2*q*c*n*b^2+2*d^2*q^2*b^2+3*d^2*q^2*b*c
*m+2*d^2*q*b*a*m*n-d^2*q*b*c*n*p+6*d^3*q*r*a-d*q^2*c^3*m-4*d*q^2*c^2*b
+2*d*q*a*b^2*m-7*d^2*q*r*c*b-d*q*c*n*b^2*m+d^2*q*b^2*n^2-3*d^4*q*r*p-2
*d*q*b^3*n-d^2*q*r*c^2*m+d^2*q^2*c^2*n-2*c^2*q^2*p*m+5*d^2*r^2*c^2-c*q
*b^3*m-q*p^3*b):}{%
}
\end{mapleinput}

\end{maplegroup}
\begin{maplegroup}
\begin{mapleinput}
\mapleinline{active}{1d}{B :=
evalf(3*m*r^2*b^3-3*r*p^2*b^3+2*m*r^3*c^2+a*d^4*q^3+d^4*r^3*m-b^3*r*n^
3+5*a*b^2*r^2-5*c^2*r*b*a^2+r^3*d^2*p-3*r*b*a^2*c*m^2+4*r*b*c*n*a^2+2*
r*p*a^2*d*m*n+d^2*q*b*a^2*m*n+d*q*a^2*b^2*m-a^2*d^2*p*c*m*q-2*d*r*a^2*
c^2*m^2+r*b*d*p*a^2-m^2*q*a^2*b^2+2*c^3*r*a^2*m+2*m^2*q*c*a^3+r*p^3*b^
2+r^2*c^3*n^2+a^2*c^2*p*d*q+a^2*c*p*b^2*m-12*a^2*c*p*n*r-3*a^2*d*q*b*c
*m^2+2*a^2*d*q*c^2*m*n-2*a^2*m^2*p^2*b*d+3*a^2*r*d^2*n^2*m-5*a^2*r*d*n
*c*m^2-10*a^2*r*b*d*m*n+4*a^2*r*b*d*m^3-a^2*d*m^2*p*b^2-2*a^2*m*p^3*c+
a^2*d^2*p*n*r+3*a^2*m*p^2*b*n+a^2*d*p*c*n*b*m+a^2*m*p^2*b*c-a^2*d*p^2*
b*n-a^2*d^2*p*b*n^2-5*a^2*d^2*p*b*q-3*a^2*d*p^2*b*c+4*a^2*c*q*b*d*n-2*
a^2*c*q*d^2*n^2+3*a^2*c^2*q*b*m-3*a^2*m*p*b^2*n+a^2*d^2*p^2*c*n-a^2*d*
p^2*c^2*m+3*a^2*d^3*p*n*q+2*a^2*b*p*c*n^2-a^2*b*p*c^2*n+2*a^2*m^2*p*c*
r+7*a^2*d^2*r*b*n-4*a^2*d^2*r*b*m^2-7*a^2*d^2*r*c*p-3*a^2*d*r*c^2*n+2*
a^2*d*p*b^2*n+6*a^2*d*r*c*n^2+3*a^2*d^3*r*m*p+5*a^2*d^2*r*n*c*m-2*a^2*
m^2*q*c^2*n+3*a^2*m^3*q*b*c+2*q*a^3*d*m*n+4*q*a^2*c*m*r+6*q*a^2*d*n*r-
8*q*a^2*m*n*r+2*q*a^2*c*n^2*d*m+q*a^2*d*p*c*m^2+5*q*a^2*d*m^2*r+2*q*c*
n*a^3+5*q*d*p*a^3-3*q*a^3*d*m^3-8*q*a^3*m*p+4*q*a^3*m^2*n+2*q*a^2*b^2*
n-q*a^3*b*m+4*q*a^2*c^2*n^2-11*q*a^2*b*r+3*a^2*d*p^3*c+2*q*a^2*b*c*p-5
*q*m^2*p*a^2*b-8*q*d^2*r*m*a^2-a^2*d*p^3*n+a^2*d*p^2*r-4*a^2*m*p^2*r+a
^2*m*p^3*d^2-a^2*b*p*n^3-4*a^2*c*q*b^2-2*a^2*c^3*q*n+a^2*m^3*p*b^2-3*a
^2*d*r*n^3-3*a^2*d^3*r*n^2-5*a^2*d*r*b^2+4*a^2*r*p*n^2-3*a^2*d^3*q^2*m
+4*a^2*d^2*q^2*c-11*a^2*d*m*r^2-2*a^2*c^2*p^2*n+3*d^3*q*r*a^2+a^2*p^4+
2*a^2*c*m*d*q^2+5*d*r*c*a^3-3*d^2*r*m*a^3+8*a^2*c^2*p*r+a^2*c*p^2*n^2+
a^2*m^2*p^2*c^2+5*a^2*m*q^2*b+2*m^3*r*a^2*c^2-2*q*a^2*c*n^3+2*q*a^2*c*
p^2-4*q*a^2*n*p^2+8*q*a^2*p*r-q*a^2*d^2*p^2+3*a^3*d*p*m^2*n+2*a^3*d*p*
c*m^2+a^5+2*a*p^2*r^2-3*a^3*d^2*p*m*n-a^3*c*n^2*d*m-d*p*a^3*b*m+6*a^2*
m*p*b*r-a^2*m^2*p*c*n*b+a^2*m*p*b*d*n^2-a^2*m*p^2*c*n*d-8*q*a^2*c*n*d*
p+3*q*a^2*d*p*n^2+q*a^2*d*p^2*m+q*a^2*n^2*b*m-8*q*a^2*c*n*b*m-q*a^2*d*
}{%
}
\end{mapleinput}

\end{maplegroup}
\begin{maplegroup}
\begin{mapleinput}
\mapleinline{active}{1d}{
m^2*n*b-3*q*a^2*d^2*p*m*n+4*q*a^2*c*n*m*p+2*q*a^2*p*b*n+4*a^3*c*n*m*p+
a^3*c*n*d*p-5*a^3*c*m*d*q-a^3*c*n*b*m+a^3*d*m^2*n*b+a^4*m^4+2*a^4*n^2+
a^3*n^4+a*q^4-a^3*m^3*n*b+8*a^3*m*n*r-a^3*d*m*n^3-3*a^3*d^2*n*q-4*a^3*
m*p*n^2-2*a^3*m^3*p*c+3*a^3*n^2*b*m+a^3*d*p*n^2-5*a^3*d*p^2*m+3*a^3*d^
2*m^2*q-2*a^3*c^2*m*p-3*q^2*a^2*c*m^2+2*q^2*c*n*a^2+3*a^4*d*m*n+m^2*p*
a^3*b+7*a^3*d*m^2*r+a^3*c*n^2*m^2-2*a^3*b*d*n^2-5*a^3*p*b*n+3*a^3*b*c*
p+4*a^3*b*d*q-6*a^3*c*m*r-7*a^3*d*n*r-q^2*d*p*a^2+4*q^2*a^2*m*p+3*a^2*
d^2*m^2*q^2-8*a^2*b*d*q^2+2*a^2*d^2*n*q^2+6*q^2*a^3-4*a^2*q^3-4*q*a^4+
5*a^3*b*r+a*d^2*n*q^3-4*a*q^2*p*r-a*d*p*q^3-a*d^3*q^3*m+10*q*d*p*a^2*b
*m-3*a^2*d^2*p*m^2*r+7*a^2*d*r*b*c*m+2*a^2*d^2*m*p^2*b-3*a*r^2*c*n^2-2
*a*d^3*r^2*p+5*a*c*r*b^3+2*a*d^2*q^2*b^2+5*a*d^2*r^2*c^2-a*q*p^3*b+5*a
*d*r^3+5*d^3*r^3*c-a*c*q^2*d^3*p+a*d*r*c^3*q-4*a*d^3*r^2*c*m-3*a*d*r*b
*c^2*m*n-9*a*d^2*r^2*c*n-2*a*d^2*r*m*c*p^2+3*a*m^2*r*b*c^2*n+4*a*d*r*p
*b*c*m^2-4*a*m^3*r*b^2*c+4*a*c*q^2*b^2+8*a*r^2*c^2*n-3*a*n*q^2*b^2+a*q
^2*c^2*n^2+3*a*d^4*r^2*n-2*a*c^4*r*p+a*b^2*q*n^3+a*r^2*d^2*q+7*a*q^2*b
*r-2*a*c^2*r*p*b*m+2*a*c*r*d^3*p^2+13*a*d^2*r^2*b*m+a*c^2*q*b*d*m*p+4*
a*c^2*r*m*p^2-2*a*c^2*r*p*d^2*n-a*b*q*c^3*p+a*b*q*d^3*p^2-a*c^2*q*p*b*
m^2+5*a*d^3*r*b*m*q+4*a*d*r*b^2*c*m^2-2*a*d^3*n*q^2*b-2*a*b^2*q*d^2*m*
p+a*d*m^2*q*b^3+2*a*d*r*c^3*m*p-2*a*b^2*q*c*n^2+3*a*b^2*q*c*d*p-6*a*c*
r*b^2*d*n+3*a*c*r*b*d^2*n^2+6*a*c*r*b*p^2-6*a*c^2*r*d*p^2+3*a*c^3*r*b*
n-5*a*m*q^2*b^2*d-8*a*m*q*b^2*r-a*q^2*c*n*d*p-4*a*d^2*q^3*c+3*a*r^2*d^
2*n^2-5*a*d^3*r^2*b-2*a*c^3*q^2*n-3*a*m*q^3*b-7*a*p*b*r^2-2*a*c*q^3*n+
4*a*q*n*r^2-4*a*m^2*r^2*c^2-a*m^3*q*b^3+a*m^2*q^2*c^3+3*a*m^2*q^2*b^2-
9*a*c*r*p*b^2+6*a*c^2*r*b*d*p-4*a*d^2*r*b*m*c*p-2*a*m^2*r*c^3*p+a*q^2*
c*p^2-5*d*r^3*c^2+a*r*d^3*q^2-2*c^4*r^2*n+5*b*c*r^3-3*b*n*r^3-a*d^2*q*
r*c^2*m+a*d^2*q^2*c^2*n-2*a*c^2*q^2*p*m-a*c*q*b^3*m+2*a*r^2*d*n*c*m+2*
a*d^2*q^2*m*b*n+3*a*m*q*b^3*n+2*a*m^2*q*b*c*r-2*a*d*q^2*b*n^2-a*m*q*b^
2*d*n^2+a*m^2*q*c*n*b^2+3*a*d^2*q^2*b*c*m-a*d^2*q*b*c*n*p-a*d*q^2*c^3*
m-4*a*d*q^2*c^2*b-7*a*d^2*q*r*c*b-a*d*q*c*n*b^2*m+a*d^2*q*b^2*n^2-3*a*
d^4*q*r*p-5*a*d*r^2*c*b-10*a*r*p*b*d*q+2*a*r*c*n*d*p^2+4*a*d^2*r^2*c*m
^2-3*a*d^3*r^2*m*n+3*a*d*r^2*c^2*m+5*a*r*b*c*n*d*p+a*r*d*q*c^2*m^2+3*a
*r*c^2*q*b+11*a*r*c*n*b^2*m+2*a*r^2*d^2*m*p-5*a*p*d*n*r^2-2*a*r*p*c^2*
n^2+5*a*b*m*n*r^2-8*a*b*d*m^2*r^2-2*a*b*c*m*r^2+4*a*c*p*d*r^2+2*a*r*c^
2*p*d*n*m-5*a*r*b*q*n^2+4*a*r*c^3*p*n+2*a*r*b*c*n*q-3*a*r*b*c*n^2*d*m+
3*a*r*b*c*n^3-6*a*r*b*c^2*n^2+2*a*r*m*q^2*c-a*q*d*m*r^2+3*a*q*d*p^2*r+
a*q*d*p*c*n*b*m+2*a*q*m*p^2*b*c-8*a*q*d*p*c*m*r-3*a*q*d^2*p*n*r+a*q*d*
p*b^2*n-3*a*q*d*p^2*b*c+a*q*d*p^2*b*n+a*b*d*n*r^2-8*a*r*p*c*n*b*m-2*a*
d*q*b^3*n-4*a*q*c^2*p*r-3*a*q*m*p*b^2*n-a*q*b*p*c*n^2+13*a*q*d*r*b^2+3
*a*c^2*p*d*q^2+2*a*q*b*p*c^2*n-a*q*d^2*r*n*c*m+11*a*q*d^2*r*c*p-2*a*q*
d*r*c^2*n+a*d^3*r*c*q*n+a*q*d*r*c*n^2+3*a*q*d^3*r*m*p-5*a*q*d^2*r*b*n-
5*a*q*d^2*r*b*m^2-a*q*c*p*b^2*m+4*a*q*c*p*n*r-a*q*d^2*m*p^2*b+2*a*q*d*
m^2*p*b^2+10*a*q*r*b*d*m*n+4*a*c*q^2*b*d*n-3*a*d*q^2*b*c*m^2-a*d*q^2*c
^2*m*n+a*c^2*q^2*b*m+a*d^2*p*b*q^2+a*q^2*c*n*b*m+a*d^2*p*c*m*q^2-a*m*p
*b*d*q^2-3*a*r*c*d*q^2-a*r*d^2*m*q^2+2*a*q*m*p*b*r+a*q^2*d*n*r+3*a*q^2
*p*b*n+3*a*c*m*d*q^3+4*a*b*d*q^3+2*d^4*r^2*b*p+4*d*r^3*b*m+5*d^2*r^2*c
*b^2+2*d*r*b^4*n-d^2*r*b^3*n^2-2*r*d*q*c*a^2+2*r*m*q*b*c^2*p-r*b*q*c*p
^2+3*r*d*p^2*b^2*c+r*b*a*n*p^2+r*b*d*p*q^2-r*b*a^2*m^2*n+r*b*a*d^2*p*m
*n-r*b*a*d*p*n^2-r*b*d^2*p*c*m*q+r*b*q*c*n*d*p-3*r*b*c^2*p*d*q-2*r*m*p
^2*b^2*c-r*d*p*b^3*n+r*a*d*m^2*n*b^2-r*a*n^2*b^2*m+2*r*m^2*p*a*b^2+r*d
^2*m*p^2*b^2+r*c*p*b^3*m-r*d*p*c*n*b^2*m-5*r*d*p*a*b^2*m-3*r*q*p*b^2*n
+r*b^2*p*c*n^2-2*r*b^2*p*c^2*n+3*r*m*p*b^3*n-r*d*p^2*b^2*n-r*d^2*p*b^2
*q+r*m*p*b^2*d*q-2*r*d*m^2*p*b^3+5*r*q*b^2*c*p+r*a*p*b^2*n-d*r^2*c^3*m
*n-4*d*r^2*b*c^2*m^2-d*r^2*p*c^2*n+r^2*c*q*d^2*n+d^2*r^2*p*c^2*m-r^2*p
*c*d*q+r^2*b*c^2*m*n+2*r^2*a*c*m*p+4*b*d*p*c*m*r^2+2*b*d^2*p*n*r^2-b*q
}{%
}
\end{mapleinput}

\end{maplegroup}
\begin{maplegroup}
\begin{mapleinput}
\mapleinline{active}{1d}{
*d*n*r^2-6*b*d^2*r^2*c*p+6*b*d*r^2*c^2*n-3*b*d*r^2*c*n^2+b*c*p*n*r^2+7
*b*r^2*c*d*q+b*r^2*d^2*m*q+3*b*d^2*r^2*n*c*m-5*m*q*b*c*r^2-d^5*r^3-7*d
*r^2*b^2*c*m+3*d^2*r^2*b^2*n+3*d*r^3*c*n-4*d^2*r^3*c*m-2*b*d*p^2*r^2+3
*b*q*p*r^2-3*b*c^2*p*r^2-d^3*q*r^2*b+7*r^2*b^2*p*d-2*r^2*q*c^2*n+3*d*r
^2*c^3*p+c^3*r^2*b*m-2*r*a*b^3*n+r*b*a^2*n^2+r*a^2*b^2*m+2*c^3*r^2*q+m
^2*r*a*b^3-3*m*r*b^4*n+2*d^2*r^2*b^2*m^2-3*m*p*b^2*r^2-5*r^2*b^2*d*m*n
+c*r^2*q^2+4*m^2*r^2*b^2*c+c*r*b^4*m-7*c*r^2*b^2*n+3*n^2*r^2*b^2+d*r*n
^2*b^3*m+d^2*r*b^2*c*n*p+d*r*c*n*b^3*m-4*d*r*b^2*c*n*q-2*d^2*r*m*q*b^2
*n+2*d*r*b^2*q*n^2+m^2*r^2*c^4-d*r*a*b^3*m+d*r*b*m*q*c^2*n-2*d^3*r^2*b
*m*p-d^2*r*b*q*c^2*n+3*d*r^2*q*m*c^2-d^3*r*b*n*p*a+4*d^2*r^2*b*c^2*m+2
*d*r*a^2*c*m*p-d^3*r^2*q*c*m-d^2*r*b*q^2*n-d^2*r*b^2*a*m*n-3*d^3*r^2*c
*b*n+2*d^3*r*n*q*b^2-d*r^2*c^4*m-5*d*r^2*c^3*b+d^2*r^2*c^3*n-c^3*r*b*q
*m^2+c^3*r*b*q*d*m+c^2*r*p*b^2*m^2+3*c*r*b^2*d*m^2*q-c^2*r*b^2*m*q-c^2
*r*b^2*d*m*p-3*c*r*b^2*d^2*m*q-d*r^3*q+m^3*r*b^4-m*r*b^2*c*n*q-m^2*r*c
*n*b^3+a*b^4*q-5*a*q^2*b*c*p-3*c*r*b*d*m*q^2-5*d*r^2*b^3-2*c*r^3*p-4*c
^2*r^2*q*d^2-d^3*r^3*n+c^5*r^2+2*a*c^2*q^3+c*r^2*d^4*q-c^2*r^2*d^3*p-2
*c^3*r^2*m*p-d*m^2*r*b^4-2*d^3*m*r^2*b^2-5*a*c^3*r^2-3*q*b^2*r^2+a*b^2
*q*c^2*n+b*r*c*q*d^3*p+2*b*r*c*q^2*n+b^2*r*c^3*p-4*b^2*r*d*q^2+5*b^3*r
*d*m*q-3*b^3*r*c*d*p+2*b^3*r*d^2*m*p+3*b^2*r*a*d^2*p+4*b^2*r*d*q*c^2-2
*b^3*r*d^2*q+3*b^3*r*n*q-3*b^3*r*m^2*q-b^3*r*c^2*n-b*r*c^4*q-b*r*d^4*q
^2-b^5*r-4*a*c^2*r*b^2*m-4*a*m*r^3+c^2*r^2*p^2+2*b*r*c^3*q*n-b*r*q*c^2
*n^2+3*r*b^4*p-2*b*r*c^2*q^2-4*b^3*r*c*q-b^2*r*d^3*p^2+b*r*d^3*q^2*m+4
*b*r*d^2*q^2*c+3*b^2*r*m*q^2+2*b^3*r*c*n^2+5*c^2*r^2*b^2-4*m*r*a^2*c^2
*n+2*q^2*a^2*n^2-4*a^2*n*r^2-a^2*d^3*p^3+5*a^2*c*r^2-4*q*a^3*n^2-4*c^2
*q^2*a^2+5*d^2*r^2*a^2-3*a^2*p^3*b+3*a^2*p^2*b^2+6*a^2*m^2*r^2+4*a^3*n
*p^2-4*a^3*p*r+a^3*d^2*n^3+3*a^3*d^2*p^2+a^3*c^2*n^2+a^4*c*m^2-a^4*b*m
-2*c*n*a^4-3*d*p*a^4-a^4*d*m^3+4*a^4*m*p-4*a^4*m^2*n+a^3*b^2*n-2*a^3*c
*n^3-4*a^3*m^3*r-3*a^3*c*p^2+a^2*c^3*p^2+2*c^2*q*a^3+2*a^3*m^2*p^2-a^2
*b^3*p+r^4-5*q^2*a^2*d*m*n+2*m*r*a^2*c*n^2+3*a*q*p^2*b^2-6*a*q*c*r^2-3
*a*q*b^3*p-2*a*c*r*p^3+a*c^4*q^2-b*r*q^3-5*d^2*r^3*b):}{%
}
\end{mapleinput}

\end{maplegroup}
\begin{maplegroup}
\begin{mapleinput}
\mapleinline{active}{1d}{s := evalf(-B/((-A)^(5/4)));}{%
}
\end{mapleinput}

\mapleresult
\begin{maplelatex}
\mapleinline{inert}{2d}{s :=
-.81584332141809188725979656230961078805698101721846912138743787689146
4709581612384523364096967240262588087176577621839593519534754857436646
38006008834788270026180162829207942917627476323600802063825066-.630576
2744559255176620934229376728312708124761783154018405570972981429024243
7736998005213095545448731240220174927694966883636580912497911035043693
968538272771370661690862684267783555502528886912001289*I;}{%
\maplemultiline{
s :=  - 
.815843321418091887259796562309610788056981017218469121387437
\backslash  \\
87689146470958161238452336409696724026258808717657762183959
\backslash  \\
35195347548574366463800600883478827002618016282920794291762
\backslash  \\
7476323600802063825066\mbox{} - 
.63057627445592551766209342293767283\backslash  \\
12708124761783154018405570972981429024243773699800521309554
\backslash  \\
54487312402201749276949668836365809124979110350436939685382
\backslash  \\
72771370661690862684267783555502528886912001289I }
}
\end{maplelatex}

\end{maplegroup}
\begin{maplegroup}
\begin{mapleinput}
\mapleinline{active}{1d}{z := -s*hypergeom([1/5,2/5,3/5,4/5],[1/2,3/4,5/4],3125/256*s^4):}{%
}
\end{mapleinput}

\end{maplegroup}
\begin{maplegroup}
\begin{mapleinput}
\mapleinline{active}{1d}{and does solve Bring's Equation}{%
}
\end{mapleinput}

\end{maplegroup}
\begin{maplegroup}
\begin{mapleinput}
\mapleinline{active}{1d}{evalf(z^5-z-s);}{%
}
\end{mapleinput}

\mapleresult
\begin{maplelatex}
\mapleinline{inert}{2d}{.2e-199*I;}{%
\[
.2\,10^{-199}\,I
\]
}
\end{maplelatex}

\end{maplegroup}
\begin{maplegroup}
\begin{mapleinput}
\mapleinline{active}{1d}{y := (-A)^(1/4)*z:}{%
}
\end{mapleinput}

\end{maplegroup}
\begin{maplegroup}
\begin{mapleinput}
\mapleinline{active}{1d}{undoing the Tschirhausian Transformation with Ferrari's method}{%
}
\end{mapleinput}

\end{maplegroup}
\begin{maplegroup}
\begin{mapleinput}
\mapleinline{active}{1d}{g :=
1/12*(-36*c*d*b-288*y*c-288*a*c+108*b^2+108*a*d^2+108*y*d^2+8*c^3+12*s
qrt(18*d^2*b^2*a+18*d^2*b^2*y+1152*d*b*y*a+240*d*b*y*c^2+240*d*b*a*c^2
-54*c*d^3*b*a-54*c*d^3*b*y-864*y*c*a*d^2+81*b^4-768*y^3-768*a^3+12*d^3
*b^3-2304*y^2*a+384*y^2*c^2-2304*y*a^2-48*y*c^4+384*a^2*c^2-48*a*c^4-3
*d^2*b^2*c^2+576*d*b*y^2+576*d*b*a^2+768*y*a*c^2-54*c*d*b^3-432*y*c*b^
2-432*y^2*c*d^2-432*a*c*b^2-432*a^2*c*d^2+162*a*d^4*y+12*a*d^2*c^3+12*
y*d^2*c^3+12*b^2*c^3+81*a^2*d^4+81*y^2*d^4))^(1/3)-12*(1/12*d*b-1/3*y-
1/3*a-1/36*c^2)/((-36*c*d*b-288*y*c-288*a*c+108*b^2+108*a*d^2+108*y*d^
2+8*c^3+12*sqrt(18*d^2*b^2*a+18*d^2*b^2*y+1152*d*b*y*a+240*d*b*y*c^2+2
40*d*b*a*c^2-54*c*d^3*b*a-54*c*d^3*b*y-864*y*c*a*d^2+81*b^4-768*y^3-76
8*a^3+12*d^3*b^3-2304*y^2*a+384*y^2*c^2-2304*y*a^2-48*y*c^4+384*a^2*c^
2-48*a*c^4-3*d^2*b^2*c^2+576*d*b*y^2+576*d*b*a^2+768*y*a*c^2-54*c*d*b^
3-432*y*c*b^2-432*y^2*c*d^2-432*a*c*b^2-432*a^2*c*d^2+162*a*d^4*y+12*a
*d^2*c^3+12*y*d^2*c^3+12*b^2*c^3+81*a^2*d^4+81*y^2*d^4))^(1/3))+1/6*c:
}{%
}
\end{mapleinput}

\end{maplegroup}
\begin{maplegroup}
\begin{mapleinput}
\mapleinline{active}{1d}{e := (d^2/4+2*g-c)^(1/2):}{%
}
\end{mapleinput}

\end{maplegroup}
\begin{maplegroup}
\begin{mapleinput}
\mapleinline{active}{1d}{f := (d*g-b)/(2*e):}{%
}
\end{mapleinput}

\end{maplegroup}
\begin{maplegroup}
\begin{mapleinput}
\mapleinline{active}{1d}{y1 := evalf(-1/4*d+1/2*e+1/4*sqrt(d^2-4*d*e+4*e^2+16*f-16*g)):}{%
}
\end{mapleinput}

\end{maplegroup}
\begin{maplegroup}
\begin{mapleinput}
\mapleinline{active}{1d}{y2 := evalf(-1/4*d+1/2*e-1/4*sqrt(d^2-4*d*e+4*e^2+16*f-16*g)):}{%
}
\end{mapleinput}

\end{maplegroup}
\begin{maplegroup}
\begin{mapleinput}
\mapleinline{active}{1d}{y3 := evalf(-1/4*d-1/2*e+1/4*sqrt(d^2+4*d*e+4*e^2-16*f-16*g)):}{%
}
\end{mapleinput}

\end{maplegroup}
\begin{maplegroup}
\begin{mapleinput}
\mapleinline{active}{1d}{y4 := evalf(-1/4*d-1/2*e-1/4*sqrt(d^2+4*d*e+4*e^2-16*f-16*g)):}{%
}
\end{mapleinput}

\end{maplegroup}
\begin{maplegroup}
\begin{mapleinput}
\mapleinline{active}{1d}{#now looking for the root that solves both the Quartic and the
Quintic}{%
}
\end{mapleinput}

\end{maplegroup}
\begin{maplegroup}
\begin{mapleinput}
\mapleinline{active}{1d}{evalf(y1^5+m*y1^4+n*y1^3+p*y1^2+q*y1+r);}{%
}
\end{mapleinput}

\mapleresult
\begin{maplelatex}
\mapleinline{inert}{2d}{-.670262027444503796935327574226396e-164-.654404080888727575792151849
0182052e-164*I;}{%
\maplemultiline{
 - .670262027444503796935327574226396\,10^{-164} \\
\mbox{} - .6544040808887275757921518490182052\,10^{-164}\,I }
}
\end{maplelatex}

\end{maplegroup}
\begin{maplegroup}
\begin{mapleinput}
\mapleinline{active}{1d}{evalf(y2^5+m*y2^4+n*y2^3+p*y2^2+q*y2+r);}{%
}
\end{mapleinput}

\mapleresult
\begin{maplelatex}
\mapleinline{inert}{2d}{782.80523423472747885830053512057667802878185155109804472138877978797
2888450235932617995871780446932630476413851455423015413001326397231266
83634661674780451644895094268241602508169015371777591202599037+3.79927
0173513100973067261361820864694346294263097462300898597218150086045413
5748601458605342688827463594103443888254180488811314627420955907635916
34218088819410407844984109417239968078437399449343508*I

;}{%
\maplemultiline{
782.80523423472747885830053512057667802878185155109804472138877978
\backslash  \\
79728884502359326179958717804469326304764138514554230154130
\backslash  \\
01326397231266836346616747804516448950942682416025081690153
\backslash  \\
71777591202599037\mbox{} + 
3.799270173513100973067261361820864694346\backslash  \\
29426309746230089859721815008604541357486014586053426888274
\backslash  \\
63594103443888254180488811314627420955907635916342180888194
\backslash  \\
10407844984109417239968078437399449343508I }
}
\end{maplelatex}

\end{maplegroup}
\begin{maplegroup}
\begin{mapleinput}
\mapleinline{active}{1d}{evalf(y3^5+m*y3^4+n*y3^3+p*y3^2+q*y3+r);}{%
}
\end{mapleinput}

\mapleresult
\begin{maplelatex}
\mapleinline{inert}{2d}{-24471.05428711853877506730011254932432187314028124956751686236513901
5142893345917771313252678577076750957332054166784797037390657774481757
164355894965388511999107977673996588972265894602403965918954074-104991
.196955226425101886129850901911773458228678168428491307753160787843223
2498327168294178778639433542591905831589669417285257948419066967216279
6753547171943376068083259295324383609618045156132*I;}{%
\maplemultiline{
 - 
24471.05428711853877506730011254932432187314028124956751686236513
\backslash  \\
90151428933459177713132526785770767509573320541667847970373
\backslash  \\
90657774481757164355894965388511999107977673996588972265894
\backslash  \\
602403965918954074\mbox{} - 
104991.196955226425101886129850901911773\backslash  \\
45822867816842849130775316078784322324983271682941787786394
\backslash  \\
33542591905831589669417285257948419066967216279675354717194
\backslash  \\
3376068083259295324383609618045156132I }
}
\end{maplelatex}

\end{maplegroup}
\begin{maplegroup}
\begin{mapleinput}
\mapleinline{active}{1d}{evalf(y4^5+m*y4^4+n*y4^3+p*y4^2+q*y4+r);}{%
}
\end{mapleinput}

\mapleresult
\begin{maplelatex}
\mapleinline{inert}{2d}{-2341.377040262803423112625582783776277016013819510961661548886093091
1512748641782916068603484330772106707019187084886849340955046992104827
929943633431723898557859216192870129460155160031024099248163764-2079.8
9649144663713312049506952871938130578787897155052839664961789785696827
6311231667458625268870521416177694775900083918775690951432710491929905
73015351337285464090257339529672681164668936088649247*I;}{%
\maplemultiline{
 - 
2341.377040262803423112625582783776277016013819510961661548886093
\backslash  \\
09115127486417829160686034843307721067070191870848868493409
\backslash  \\
55046992104827929943633431723898557859216192870129460155160
\backslash  \\
031024099248163764\mbox{} - 
2079.89649144663713312049506952871938130\backslash  \\
57878789715505283966496178978569682763112316674586252688705
\backslash  \\
21416177694775900083918775690951432710491929905730153513372
\backslash  \\
85464090257339529672681164668936088649247I }
}
\end{maplelatex}

\end{maplegroup}
\begin{maplegroup}
\begin{mapleinput}
\mapleinline{active}{1d}{in this example y1 is the root we want, let it be the first root of
the Quintic, r1}{%
}
\end{mapleinput}

\end{maplegroup}
\begin{maplegroup}
\begin{mapleinput}
\mapleinline{active}{1d}{r1 := y1:}{%
}
\end{mapleinput}

\end{maplegroup}
\begin{maplegroup}
\begin{mapleinput}
\mapleinline{active}{1d}{factoring it out of the Quintic, leaving only the 
Quartic to solve}{%
}
\end{mapleinput}

\end{maplegroup}
\begin{maplegroup}
\begin{mapleinput}
\mapleinline{active}{1d}{dd := m+r1:}{%
}
\end{mapleinput}

\end{maplegroup}
\begin{maplegroup}
\begin{mapleinput}
\mapleinline{active}{1d}{cc := n+r1^2+m*r1:}{%
}
\end{mapleinput}

\end{maplegroup}
\begin{maplegroup}
\begin{mapleinput}
\mapleinline{active}{1d}{bb := p+r1*n+r1^3+m*r1^2:}{%
}
\end{mapleinput}

\end{maplegroup}
\begin{maplegroup}
\begin{mapleinput}
\mapleinline{active}{1d}{aa := q+r1*p+r1^2*n+r1^4+m*r1^3:}{%
}
\end{mapleinput}

\end{maplegroup}
\begin{maplegroup}
\begin{mapleinput}
\mapleinline{active}{1d}{yy := 0:}{%
}
\end{mapleinput}

\end{maplegroup}
\begin{maplegroup}
\begin{mapleinput}
\mapleinline{active}{1d}{gg :=
1/12*(-36*cc*dd*bb-288*yy*cc-288*aa*cc+108*bb^2+108*aa*dd^2+108*yy*dd^
2+8*cc^3+12*sqrt(18*dd^2*bb^2*aa+18*dd^2*bb^2*yy+1152*dd*bb*yy*aa+240*
dd*bb*yy*cc^2+240*dd*bb*aa*cc^2-54*cc*dd^3*bb*aa-54*cc*dd^3*bb*yy-864*
yy*cc*aa*dd^2+81*bb^4-768*yy^3-768*aa^3+12*dd^3*bb^3-2304*yy^2*aa+384*
yy^2*cc^2-2304*yy*aa^2-48*yy*cc^4+384*aa^2*cc^2-48*aa*cc^4-3*dd^2*bb^2
*cc^2+576*dd*bb*yy^2+576*dd*bb*aa^2+768*yy*aa*cc^2-54*cc*dd*bb^3-432*y
y*cc*bb^2-432*yy^2*cc*dd^2-432*aa*cc*bb^2-432*aa^2*cc*dd^2+162*aa*dd^4
*yy+12*aa*dd^2*cc^3+12*yy*dd^2*cc^3+12*bb^2*cc^3+81*aa^2*dd^4+81*yy^2*
dd^4))^(1/3)-12*(1/12*dd*bb-1/3*yy-1/3*aa-1/36*cc^2)/((-36*cc*dd*bb-28
8*yy*cc-288*aa*cc+108*bb^2+108*aa*dd^2+108*yy*dd^2+8*cc^3+12*sqrt(18*d
d^2*bb^2*aa+18*dd^2*bb^2*yy+1152*dd*bb*yy*aa+240*dd*bb*yy*cc^2+240*dd*
bb*aa*cc^2-54*cc*dd^3*bb*aa-54*cc*dd^3*bb*yy-864*yy*cc*aa*dd^2+81*bb^4
-768*yy^3-768*aa^3+12*dd^3*bb^3-2304*yy^2*aa+384*yy^2*cc^2-2304*yy*aa^
2-48*yy*cc^4+384*aa^2*cc^2-48*aa*cc^4-3*dd^2*bb^2*cc^2+576*dd*bb*yy^2+
576*dd*bb*aa^2+768*yy*aa*cc^2-54*cc*dd*bb^3-432*yy*cc*bb^2-432*yy^2*cc
*dd^2-432*aa*cc*bb^2-432*aa^2*cc*dd^2+162*aa*dd^4*yy+12*aa*dd^2*cc^3+1
2*yy*dd^2*cc^3+12*bb^2*cc^3+81*aa^2*dd^4+81*yy^2*dd^4))^(1/3))+1/6*cc:
}{%
}
\end{mapleinput}

\end{maplegroup}
\begin{maplegroup}
\begin{mapleinput}
\mapleinline{active}{1d}{ee := (dd^2/4+2*gg-cc)^(1/2):}{%
}
\end{mapleinput}

\end{maplegroup}
\begin{maplegroup}
\begin{mapleinput}
\mapleinline{active}{1d}{ff := (dd*gg-bb)/(2*ee):}{%
}
\end{mapleinput}

\end{maplegroup}
\begin{maplegroup}
\begin{mapleinput}
\mapleinline{active}{1d}{yy1 :=
evalf(-1/4*dd+1/2*ee+1/4*sqrt(dd^2-4*dd*ee+4*ee^2+16*ff-16*gg)):}{%
}
\end{mapleinput}

\end{maplegroup}
\begin{maplegroup}
\begin{mapleinput}
\mapleinline{active}{1d}{yy2
:=evalf(-1/4*dd+1/2*ee-1/4*sqrt(dd^2-4*dd*ee+4*ee^2+16*ff-16*gg)):}{%
}
\end{mapleinput}

\end{maplegroup}
\begin{maplegroup}
\begin{mapleinput}
\mapleinline{active}{1d}{yy3 :=evalf(
-1/4*dd-1/2*ee+1/4*sqrt(dd^2+4*dd*ee+4*ee^2-16*ff-16*gg)):}{%
}
\end{mapleinput}

\end{maplegroup}
\begin{maplegroup}
\begin{mapleinput}
\mapleinline{active}{1d}{yy4
:=evalf(-1/4*dd-1/2*ee-1/4*sqrt(dd^2+4*dd*ee+4*ee^2-16*ff-16*gg)):}{%
}
\end{mapleinput}

\end{maplegroup}
\begin{maplegroup}
\begin{mapleinput}
\mapleinline{active}{1d}{Do the roots of the Quartic satisfy the Quintic?}{%
}
\end{mapleinput}

\end{maplegroup}
\begin{maplegroup}
\begin{mapleinput}
\mapleinline{active}{1d}{evalf(yy1^5+m*yy1^4+n*yy1^3+p*yy1^2+q*yy1+r);}{%
}
\end{mapleinput}

\mapleresult
\begin{maplelatex}
\mapleinline{inert}{2d}{-.670262027444503796935327574226455e-164-.654404080888727575792151848
993e-164*I;}{%
\maplemultiline{
 - .670262027444503796935327574226455\,10^{-164} \\
\mbox{} - .654404080888727575792151848993\,10^{-164}\,I }
}
\end{maplelatex}

\end{maplegroup}
\begin{maplegroup}
\begin{mapleinput}
\mapleinline{active}{1d}{evalf(yy2^5+m*yy2^4+n*yy2^3+p*yy2^2+q*yy2+r);}{%
}
\end{mapleinput}

\mapleresult
\begin{maplelatex}
\mapleinline{inert}{2d}{-.670262027444503796935327575500680e-164-.654404080888727575792151e-1
64*I;}{%
\maplemultiline{
 - .670262027444503796935327575500680\,10^{-164} \\
\mbox{} - .654404080888727575792151\,10^{-164}\,I }
}
\end{maplelatex}

\end{maplegroup}
\begin{maplegroup}
\begin{mapleinput}
\mapleinline{active}{1d}{evalf(yy3^5+m*yy3^4+n*yy3^3+p*yy3^2+q*yy3+r);}{%
}
\end{mapleinput}

\mapleresult
\begin{maplelatex}
\mapleinline{inert}{2d}{-.670262027444503796935327574225863e-164-.654404080888727575792151848
9997934e-164*I;}{%
\maplemultiline{
 - .670262027444503796935327574225863\,10^{-164} \\
\mbox{} - .6544040808887275757921518489997934\,10^{-164}\,I }
}
\end{maplelatex}

\end{maplegroup}
\begin{maplegroup}
\begin{mapleinput}
\mapleinline{active}{1d}{evalf(yy4^5+m*yy4^4+n*yy4^3+p*yy4^2+q*yy4+r);}{%
}
\end{mapleinput}

\mapleresult
\begin{maplelatex}
\mapleinline{inert}{2d}{-.670262027444503796935327574247721e-164-.654404080888727575792151849
0115065e-164*I;}{%
\maplemultiline{
 - .670262027444503796935327574247721\,10^{-164} \\
\mbox{} - .6544040808887275757921518490115065\,10^{-164}\,I }
}
\end{maplelatex}

\end{maplegroup}
\begin{maplegroup}
\begin{mapleinput}
\mapleinline{active}{1d}{They do. The five roots of the General Quintic equation are;}{%
}
\end{mapleinput}

\end{maplegroup}
\begin{maplegroup}
\begin{mapleinput}
\mapleinline{active}{1d}{r1 := r1;}{%
}
\end{mapleinput}

\mapleresult
\begin{maplelatex}
\mapleinline{inert}{2d}{r1 :=
.377479231046779905370547206961209593544489801159674781330119420080378
8972241762228066524401453869059411776882337250030873958594005806297748
3352794710876598360361485361542403140749630830436441090570208-.4668210
0370132134176098284066265750022049699904489626603055650351534241700234
0501534541218808904885455382112417283951721517559653844115256905396985
732092526163375444738646870686021501830132473555311*I;}{%
\maplemultiline{
\mathit{r1} := 
.3774792310467799053705472069612095935444898011596747813301194
\backslash  \\
20080378897224176222806652440145386905941177688233725003087
\backslash  \\
39585940058062977483352794710876598360361485361542403140749
\backslash  \\
630830436441090570208\mbox{} - 
.466821003701321341760982840662657500\backslash  \\
22049699904489626603055650351534241700234050153454121880890
\backslash  \\
48854553821124172839517215175596538441152569053969857320925
\backslash  \\
26163375444738646870686021501830132473555311I }
}
\end{maplelatex}

\end{maplegroup}
\begin{maplegroup}
\begin{mapleinput}
\mapleinline{active}{1d}{r2 := yy1;}{%
}
\end{mapleinput}

\mapleresult
\begin{maplelatex}
\mapleinline{inert}{2d}{r2 :=
.311526239599294542035310266802123056675813293069289355469523096066004
5088256528877398989517046008082026041181960272271663366901523137479556
0798035158860462565687955370738966390212760479641932886238545e-2-6.515
8800625716849604275554567157395619099490178856720500232969490753571089
8235149830470703576867043470529449191796836836837700865471492870192017
5749825825560080670253184158055421829508023205091953137*I;}{%
\maplemultiline{
\mathit{r2} := 
.0031152623959929454203531026680212305667581329306928935546952
\backslash  \\
30960660045088256528877398989517046008082026041181960272271
\backslash  \\
66336690152313747955607980351588604625656879553707389663902
\backslash  \\
12760479641932886238545\mbox{} - 
6.515880062571684960427555456715739\backslash  \\
56190994901788567205002329694907535710898235149830470703576
\backslash  \\
86704347052944919179683683683770086547149287019201757498258
\backslash  \\
25560080670253184158055421829508023205091953137I }
}
\end{maplelatex}

\end{maplegroup}
\begin{maplegroup}
\begin{mapleinput}
\mapleinline{active}{1d}{r3 := yy2;}{%
}
\end{mapleinput}

\mapleresult
\begin{maplelatex}
\mapleinline{inert}{2d}{r3 :=
.280621403507492946823494953041804575855110247063853745991415245436122
7402630026122602677255801204011950151918487138814436726913624853479872
139207968095758568959370275843832797806302776026675927402455e-3+206.48
9462435352478709772156013095764283126100332544514609464387704273436398
8734937027558150313243206129966052604378809617037327210535844703819997
6185463002277378139422837250633105266650324928524956067*I;}{%
\maplemultiline{
\mathit{r3} := 
.0002806214035074929468234949530418045758551102470638537459914
\backslash  \\
15245436122740263002612260267725580120401195015191848713881
\backslash  \\
44367269136248534798721392079680957585689593702758438327978
\backslash  \\
06302776026675927402455\mbox{} + 
206.4894624353524787097721560130957\backslash  \\
64283126100332544514609464387704273436398873493702755815031
\backslash  \\
32432061299660526043788096170373272105358447038199976185463
\backslash  \\
002277378139422837250633105266650324928524956067I }
}
\end{maplelatex}

\end{maplegroup}
\begin{maplegroup}
\begin{mapleinput}
\mapleinline{active}{1d}{r4 := yy3;}{%
}
\end{mapleinput}

\mapleresult
\begin{maplelatex}
\mapleinline{inert}{2d}{r4 :=
.345083782925994921623416729927043173122111832720745591758875859061234
8400205304806880973384062988220741607533348718973543601293057300924385
6080216228866292791120138190327894714179764114732975550346486+.4636972
2578119903036077263312428888132284271256780564831601970143334073300130
8691601930073641175868971448605695998561352418556148809754449142410912
841189847519928621430651231210260275880679737044336*I;}{%
\maplemultiline{
\mathit{r4} := 
.3450837829259949216234167299270431731221118327207455917588758
\backslash  \\
59061234840020530480688097338406298822074160753334871897354
\backslash  \\
36012930573009243856080216228866292791120138190327894714179
\backslash  \\
764114732975550346486\mbox{} + 
.463697225781199030360772633124288881\backslash  \\
32284271256780564831601970143334073300130869160193007364117
\backslash  \\
58689714486056959985613524185561488097544491424109128411898
\backslash  \\
47519928621430651231210260275880679737044336I }
}
\end{maplelatex}

\end{maplegroup}
\begin{maplegroup}
\begin{mapleinput}
\mapleinline{active}{1d}{r5 := yy4;}{%
}
\end{mapleinput}

\mapleresult
\begin{maplelatex}
\mapleinline{inert}{2d}{r5 :=
-.72595889777227526536114053450931580180921487705817712038968192534770
9905073226234984409035794311856498559497942405886594863028299196345040
93762383371012453362828096808336125846809585577726102729053104+.295414
0513932856205560965115834389768150297181824805827344604688392239410988
9605481503149612078838192622738201372797033743291106164819127683399968
086705464942722848119647363844924230723372578903450e-1*I;}{%
\maplemultiline{
\mathit{r5} :=  - 
.72595889777227526536114053450931580180921487705817712038968
\backslash  \\
19253477099050732262349844090357943118564985594979424058865
\backslash  \\
94863028299196345040937623833710124533628280968083361258468
\backslash  \\
09585577726102729053104\mbox{} + 
.0295414051393285620556096511583438\backslash  \\
97681502971818248058273446046883922394109889605481503149612
\backslash  \\
07883819262273820137279703374329110616481912768339996808670
\backslash  \\
5464942722848119647363844924230723372578903450I }
}
\end{maplelatex}

\end{maplegroup}

\begin{maplegroup}

\end{maplegroup}
\vspace{0.3cm}
Acknowlegement: Andrew DeBenedictis for helping in the latex preparation
of this document.

\begin{maplegroup}
\begin{Heading 1}
REFERENCES
\end{Heading 1}

\end{maplegroup}
\begin{maplegroup}
1. Young, Robyn V. "Notable Mathematicians",Gale, Detroit, 1998

\end{maplegroup}
\begin{maplegroup}
2. Guerlac, Henry."Biographical Dictionary of Mathematics", Collier    
       Macmillan Canada, Vol. I to IV, 1991.

\end{maplegroup}
\begin{maplegroup}
3. Bring, Erland Sam."B. cum D. Meletemata quae dam mathematica circa
\#        transformation aequationum algebraicarum, quae consent".
Ampliss. Facult. Philos. in Regia Academia Carolina Praeside D. Erland
Sam. Bring, Hist. Profess-Reg. \& Ord.publico Eruditorum Examini
modeste subjicit Sven Gustaf Sommelius, Stipendiarius Regius \&
Palmcreutzianus Lundensis.Die XIV. The main part of Bring's work is
reproduced in Quar. Jour. Math., 6, 1864; Archiv. Math.
Phys., 41, 1864, 105-112; Annali di Mat., 6, 1864, 33-42. There also a
republication of his works in 12 volumes at the University of Lund.

\end{maplegroup}
\begin{maplegroup}
4. Harley, Rev. Robert, "A Contribution to the History of the Problem
of the Reduction of the General Equation of the Fifth Degree to a
Trinomial Form". Quar. Jour. Math., 6, 1864, pp. 38-47.

\end{maplegroup}
\begin{maplegroup}
5. Weisstein, Eric W. "CRC Concise Encyclopedia of Mathematics", CRC
Press, 1999. pp. 1497 to 1500.

\end{maplegroup}
\begin{maplegroup}
6. King, Bruce R. "Beyond the Quartic Equation", Birkhauser, Boston,
1996, pp. 87-89.

\end{maplegroup}
\begin{maplegroup}
7. Klein, Felix. "Lectures on the Icosahedron and the Solution of
Equations of the Fifth Degree" Dover Publishing, NY.NY, 1956.

\end{maplegroup}
\begin{maplegroup}
8. Prasolov, Viktor and Solovyev,"Elliptic Functions and Elliptic
Integrals", Translation of Mathematical Monographs, Vol. 170. American
Mathematic Society, Providence, Rhode Island, 1997.

\end{maplegroup}
\begin{maplegroup}
9. Cockle, James. "Sketch of a Theory of Transcendental Roots". Phil.
Mag. Vol 20, pp. 145-148,1860.

\end{maplegroup}
\begin{maplegroup}
10. Cockle, James, "On Transcendental and Algebraic
Solution.-Supplementary Paper". Phil. Mag. Vol. 13, pp. 135-139.

\end{maplegroup}
\begin{maplegroup}
11. Harley, Rev Robert. "On the theory of the Transcendental Solution
of Algebraic Equations". Quart. Journal of Pure and Applied Math, Vol.
5 p. 337. 1862.

\end{maplegroup}
\begin{maplegroup}
12. Cayley, Arthur. "Note on a Differential Equation", Memoirs of the
Literary and Philosophical Society of Manchester, vol. II (1865), pp.
111-114. Read February 18,1862.

\end{maplegroup}
\begin{maplegroup}
13. Cayley, Arthur. "On Tschirnhaus's Tranformation". Phil. Trans.
Roy. Soc. London, Vol. 151. pp. 561-578. Also in Collected Mathematical
Papers, Cayley p. 275. 

\end{maplegroup}
\begin{maplegroup}
14. Green, Mark L. "On the Analytic Solution of the Equation of Fifth
Degree", Compositio Mathematica, Vol. 37, Fasc. 3, 1978. p. 233-241.
Sijthoff \& Noordhoff International Publishers-Alphen aan den Rijn.
Netherlands.

\end{maplegroup}
\begin{maplegroup}
15. Slater, Lucy Joan. "Generalized Hypergeometric Functions",
Cambridge University Press, 1966, pp. 42-44.

\end{maplegroup}
\end{document}